\documentclass[12pt]{article}
\usepackage{amsmath,amsfonts,amssymb,amsthm}

\newcommand{\I}{{\bf 1}}
\newtheorem{proposition}{Proposition}[section]
\newtheorem{theorem}[proposition]{Theorem}

\newtheorem{lemma}[proposition]{Lemma}

\newcommand{\nc}{\newcommand}
\nc{\R}{{\mathbb R}}
\nc{\N}{{\mathbb N}}
\nc{\Z}{{\mathbb Z}}
\DeclareMathOperator{\card}{card}
\nc{\BP}{\mathbb{P}}
\nc{\BE}{\mathbb{E}}
\nc{\BQ}{\mathbb{Q}}
\nc{\bN}{{\mathbf N}}
\nc{\bU}{{\mathbf U}}

\newcommand{\eqco}{\setcounter{equation}{0}}

\newcommand{\eps}{\varepsilon}
\newcommand{\bea}{\begin{eqnarray}}
\newcommand{\eea}{\end{eqnarray}}
\newcommand{\bean}{\begin{eqnarray*}}
\newcommand{\eean}{\end{eqnarray*}}
\newcommand{\diam}{{\rm diam}}
\newcommand{\X}{\chi}
\newcommand{\Rex}{R_{{\rm ex}}}
\newcommand{\Var}{{\rm Var}}
\newcommand{\tod}{\stackrel{{\cal D}}{\longrightarrow}}
\newcommand{\G}{{\cal G}}

\newcommand{\lbl}{\label}

\newcommand{\eq}[1]{$(\ref{#1})$}

\begin{document}

\author{G\"unter Last\footnote{
Institut f\"ur Stochastik,
 Karlsruher Institut f\"ur Technologie,
76128 Karlsruhe, Germany. 
Email: guenter.last@kit.edu}
\ and 
\ Mathew D. Penrose\footnote{ Department of Mathematical Sciences, University of Bath,
Bath BA2 7AY, United Kingdom,
Email: m.d.penrose@bath.ac.uk} \footnote{Partially
 supported by the Alexander von Humboldt Foundation
through a Friedrich Wilhelm Bessel Research Award, and by the Isaac Newton
Institute for Mathematical Sciences, Cambridge}
}
\title{Percolation and limit theory for the Poisson lilypond model}
\date{\today}
\maketitle
\begin{abstract}
\noindent
The lilypond model on a point process in $d$-space is a growth-maximal
system of non-overlapping balls centred at the points. We establish
central limit theorems for the total volume and the number of components
of the lilypond model on a sequence of Poisson or binomial point 
processes on expanding windows. For the lilypond
model over a homogeneous Poisson process, we  give subexponentially
decaying tail bounds for the size
of the cluster at the origin. Finally, we consider
the enhanced Poisson lilypond model where  all the balls are 
enlarged by a fixed amount (the enhancement parameter),   
and show that for $d > 1 $ the 
critical value of this parameter,
above which the enhanced model percolates, 
is strictly positive.

\end{abstract}

\section{Introduction}
\lbl{secintro}
\setcounter{equation}{0}
Suppose $\varphi$ is a locally finite set of points
of cardinality at least $2$ in the space $\R^d$.
The {\em lilypond model} based on $\varphi$ is
the system of balls (or {\em grains})
 $\{B_{\rho(x)}(x):x\in\varphi\}$
(here $B_r(x): = \{y \in \R^d: |y-x| \leq r\}$ and $|\cdot|$ is
Euclidean norm)
with the following two properties:
\begin{itemize}
\item The {\em hard-core property}: 
 $\rho(x) + \rho(y) \leq |x-y|$ 
for all different $x,y\in\varphi$. 
\item The {\em smaller grain-neighbour property}:
for each $x\in\varphi$ there is at least one
$y\in\varphi\setminus\{x\}$ such that
$\rho(x) + \rho(y) = |y-x| $
(in which case the points $x$ and $y$ are called {\em grain-neighbours})
and $\rho(y)\le \rho(x)$. In this case we call $y$ a
{\em smaller} grain-neighbour of $x$.
\end{itemize}
In the sequel we shall write $\rho(x,\varphi)$ to denote the
dependence of the radii on both $\varphi$ and $x\in \varphi$.
Heveling and Last \cite{HL06} established 
 existence and uniqueness
of the model for all such $\varphi$ (in fact in greater
generality). Define the union set
\begin{align}
Z(\varphi) :=\bigcup_{x\in \varphi} B_{\rho(x,\varphi)}(x).
\end{align}
In the case of finite $\varphi$, the lilypond model may be
 constructed as follows. All points of $\varphi$
start growing at the same time and at the same rate.
Any given ball ceases its growth as soon as it reaches any other ball.
When $\varphi$ has just a single element $x$, we define
$\rho(x,\varphi):= +\infty$.

The model is of interest since it is a {\em growth-maximal hard-core}
model. No single ball can grow further without
overlapping another ball.  
It  also has a {\em maximin} property: 
if  $\varphi$ is finite
with $n$ elements and the radii $\{\rho(x,\varphi): x \in \varphi\}$
are listed in ascending order as $\rho_1,\rho_2, \ldots \rho_n$,
and if
$\{\rho'(x):x \in \varphi\}$ is any other system of radii 
satisfying the hard-core property, similarly listed
 in ascending order as   $\rho'_1, \ldots,\rho'_{n}$
then $(\rho_1,\ldots,\rho_n)$ exceeds $(\rho'_1,\ldots,\rho'_n)$
in the lexicographic ordering on $\R^n$.

In this paper we consider the lilypond model on
{\em random} $\varphi$.
The lilypond model was introduced
by  H\"aggstr\"om and Meester \cite{HM96} 
for the case where $\varphi$ is a stationary
Poisson process $\Phi$ in $\R^d$ of intensity one; 
 we call this the {\em Poisson lilypond model} and set
$\Psi=\{(x,\rho(x,\Phi)):x\in\Phi\}$.
They proved that the union set $Z:=Z(\Phi)$
does not percolate, i.e.\ does not have an unbounded connected
component. Interestingly, there {\em does} exist a
stationary percolating hard-core system of (non-lilypond) grains on $\Phi$, at 
least in high dimensions; see \cite{CHR}.

Apart from the one-dimensional case (see Daley, Mallows and Shepp \cite{DMS00}),
only a few further probabilistic properties of
$\Psi$ are known.
 Daley, Stoyan and Stoyan \cite{DSS99} give bounds
for the {\em volume fraction}
\begin{align}
p_Z:=\BP(0\in Z)=\BE V_d(Z\cap[0,1]^d)
\end{align}
of $Z$, where $0$ denotes the origin in $\R^d$ and
$V_d$ is Lebesgue measure (volume) on $\R^d$.
The latter paper  also has some numerical results
on the {\em typical radius} $\rho_0$ of
$\Psi$. The distribution of $\rho_0$ is that of $\rho(x,\Phi)$
for a ``randomly picked'' $x\in\Phi$.
Because $\Phi$ is Poisson, it is well-known that
the distribution of $\rho_0$ is that of
$\rho(0,\Phi^0)$ where $\Phi^0:=\Phi\cup\{0\}$. 
Because of the hard-core property of $\Psi$ we clearly
have
\begin{align}\label{pZ}
p_Z=b_d\BE[\rho_0^d],
\end{align}
where $b_d:=V_d(B_1(0))$ is the volume of the unit ball in $\R^d$.

The contributions of the present paper fall into three categories:
tail bounds, central limit theorems, and non-percolation
under positive enhancement.
These may be viewed as extending
the percolation theory of the Poisson lilypond model beyond the basic
fact that $Z$ does not percolate. We now give an overview
of our results, which are proved using notions
of {\em stabilization} developed in Sections \ref{sstable}
and \ref{secstrong}.

Given $\varphi$ as above,
and given $x \in \varphi$,
define $C'(x,\varphi)$ to be the connected component
of $Z(\varphi)$ containing $x$.  We study the `cluster at the origin'
 $C'(0,\Phi^0)$, which amounts to studying
 $C'(x,\Phi)$ for a randomly picked point $x\in \Phi$. Since
 $Z$ does not percolate, we know that $C'(0,\Phi^0)$ 
is almost surely bounded.  In Section \ref{secdecay}
we shall give tail bounds on its size.
We consider three different measures of the `size' of $C'(0,\Phi^0)$,
 namely volume, metric diameter,
and the number of constituent grains.
In each case, we
shall give upper tail bounds showing that the probability
of the size of a cluster exceeding $t$ decays exponentially
in a power of $t$. We also give lower bounds of the same form
but with different exponents; it remains open to 
establish the `correct' exponent (if any) for the tail decay 
of the size of the cluster at the origin.   

Our central limit theorems are stated with reference to a sequence
of expanding windows in $\R^d$.
Let $W\subset\R^d$ be a compact convex set  containing an open neighborhood
of the origin with $V_d(W) =1$, and set $W_n := n^{1/d} W$,
where for $A \subset \R^d$ and $x \in \R^d$, $a \in \R_+$ 
we write $x + aA $ for $\{x+ay: y \in A\}$.
For $ A \subset \R^d$,
let $\kappa(A)$ denote the number of connected
components of $A$.  In Section \ref{secCLT}, we shall
 derive  central limit theorems for (among other things) $V_d(Z(\Phi_n))$,
and $\kappa(Z(\Phi_n))$,
where we set  $\Phi_n := \Phi \cap W_n$.
We shall also establish the corresponding de-Poissonized central limit
theorems, where instead of $\Phi_n$ one considers 
a point process $\chi_{n}$ consisting of $n$ independent
uniformly distributed  points in $n^{1/d} W$.
Central limit theorems such as these could be of use in
establishing confidence intervals for quantities such as the volume
fraction, based on simulations.

In Section \ref{secfrog}
we consider  the {\em enhanced lilypond model}
 $Z^\delta$, defined for $\delta>0$
to be the $\delta$-neighbourhood of $Z$ in $\R^d$. We shall
show that for $d \geq 2$, there is a strictly positive critical value of
$\delta$ such that the enhancement   $Z^\delta$
percolates almost surely if $\delta$ is above this critical value,
and does not percolate if $\delta$ is below the critical value.
This gives us a one-parameter family of random sets in
the continuum exhibiting a non-trivial phase transition. 

As a final remark here, we compare the Poisson
lilypond model
to  the {\em random sequential adsorption} (RSA) process with
infinite input, in which unit balls 
arrive at locations given by a homogeneous space-time
Poisson process starting at time zero, each ball
being irreversibly accepted if it does not overlap any previously accepted
ball. Continuing to time infinity, one ends up with 
a maximal system of balls satisfying the hard-core 
constraint, as with the lilypond model; see Penrose
\cite{PeRSA}, and Schreiber, Penrose and Yukich
\cite{SPY} for formal definitions, existence
and limit theorems for RSA.  Noteworthy
 differences are that for RSA, the radii are all the same, 
 the point process of accepted ball centres is not 
a spatial Poisson process, and the balls almost surely  do not touch.

Notation: 
 We use $c, c', c''$ and so on to denote positive
finite constants whose values are unimportant and may change from line to line.
On  the other hand, we denote by $c_1,c_2$ and so on, constants
whose values (though still not very important)
are fixed and which may reappear in other parts of the paper. 
 For nonempty $A \subset \R^d$ we write $\diam(A)$
for $\sup_{(x,y) \in A \times  A} |x-y|$.
We write $\card(A)$ for the number of elements of $A$ (possibly $\infty$).

\section{Stabilization}\label{sstable}
\setcounter{equation}{0}

In this section, we establish that
 there is an almost surely finite random variable $R := R(\Phi)$
such that the radius $\rho(0,\Phi \cup \{0\})$ is unaffected by modifications
to the point set $\Phi$ outside the ball $B_R(0)$.
This is known as a {\em radius of stabilization} for $\rho(0,\Phi)$.
Moreover, we establish tail bounds for $R$, i.e. bounds
 on the $\BP(R>t)$, which decay exponentially in $t^{d/(d+1)}$.
%
%
These tail
estimates will be crucial in all of our subsequent proofs. 

A {\em point process} is defined as
a random variable taking values in the space $\bN$ of all locally
finite subsets of $\R^d$ equipped with the smallest $\sigma$-field
$\mathcal{N}$ containing the sets $\{\varphi\in\bN:\varphi(B)=k\}$
for all Borel $B \subset  \R^d$
 and all $k\in\N_0$, where
$\varphi(B)$ denotes the number of elements of $\varphi \cap B$. 
By a (finite) {\em descending chain} in $\varphi\in\bN$ we mean
a finite sequence $x_0,\ldots,x_n$ ($n\ge 1$) of distinct points
of $\varphi$ for which $|x_{i-1}-x_i|\ge |x_i-x_{i+1}|$ for all 
$i\in\{1,\ldots,n-1\}$.
Note that any two points of $\varphi$, considered 
on their own, form a descending
chain. Let $\varphi\in\bN$
and $x\in\varphi$. If $\varphi(\R^d) \geq 2$
 we let $D(x,\varphi)$ denote the
{\em nearest neighbour distance}
from $x$ in $\varphi\setminus\{x\}$. 
That is, we set
$D(x,\varphi) : =\min\{|x-y|:y\in \varphi\setminus\{x\}\}$.
If $\varphi =\{x\}$, set $D(x,\varphi):= +\infty$.

For any $\varphi\in\bN$ and $x,y  \in \R^d$
we set $\varphi^x:=\varphi\cup\{x\}$ and $\varphi^{x,y} = 
\varphi \cup \{x,y\}$.
We construct a closed set $S(y,\varphi) \subset \R^d$ such that 
if this set is bounded, then
the radius $\rho(y,\varphi^y)$
is determined by the restriction of $\varphi$ to $S(y,\varphi)$ 
(see \eq{2.1} below).
In the trivial case where  $\varphi \setminus \{y\} = \emptyset$
 we define 
$S(y,\varphi):=\R^d$. 
Otherwise, we define
\begin{align}\label{2.2}
S(y,\varphi):= 
B_{2D(y,\varphi^y)}(y)
 \cup \bigcup_{(x,r)\in A(y,\varphi) }B_{r}(x),
\end{align}
where the set $A(y,\varphi)\subset\varphi\times(0,\infty)$
is defined as follows. A pair $(x,r)$ belongs to $A(\varphi)$
if there is a descending chain $x_0,\ldots,x_n$ in $\varphi^0$
such that $x_0=y$, $|x_1-y|\le 2D(y,\varphi^y)$,
$x_n=x$ and $r=|x-x_{n-1}|$.

The next result tells us essentially that $S(y,\Phi)$ is a {\em stopping set};
see  \cite{BaLa09}, \cite{Zuyev}.

\begin{lemma}
\lbl{stopsetlem}
Suppose $y \in \R^d$ and $\varphi, \psi  \in \bN$ with
 $\psi \cap S(y,\varphi) = \varphi \cap S(y,\varphi)$.
Then $S(y,\psi)= S(y,\varphi)$.
\end{lemma}

{\em Proof.} 
Assume  $\varphi \setminus \{y\}$ is nonempty (otherwise the
result is trivial).  Then $D(y,\varphi^y) = D(y,\psi^y)$.

Suppose $(x,s) \in A(y,\varphi)$.  Then
there is descending chain $x_0,\ldots,x_n$ in $\varphi^y$
such that $x_0=y$, $|x_1-y|\le 2D(y,\varphi^y)$,
$x_n=x$ and $s=|x-x_{n-1}|$.
For $1 \leq m \leq n$ we have $x_m \in S(y,\varphi) \cap \varphi$
 so that
$x_m \in \psi$.  
%
%
Hence, 
 $(x,s) \in  A(y,\psi)$,
so $A(y,\varphi) \subset A(y,\psi)$.

Conversely, suppose $(x,s) \in A(y,\psi)$. 
Then there is descending chain $x_0,\ldots,x_n$ in $\psi^y$
such that $x_0=y$, $|x_1-y|\le 2D(y,\psi^y)$,
$x_n=x$ and $s=|x-x_{n-1}|$.

We claim that
 $x_m \in \varphi$ 
for $1 \leq m \leq n$.
This is proved by induction on $m$; if it holds for $1 \leq m \leq k-1$ 
then $x_k \in S(y,\varphi) \cap \psi$ so $x_k \in \varphi$. 
To start the induction note that $|x_1-y| \leq 2 D(y,\varphi^y)$
so $x_1 \in S(y,\varphi) \cap \psi$ so $x_1 \in \varphi$.

By the preceding claim,
 $(x,s) \in  A(y,\varphi)$ 
so $A(y,\psi) \subset A(y,\varphi)$, and hence 
$S(y,\psi) \subset S(y,\varphi)$.
Therefore $A(y,\psi) = A(y,\varphi)$, so that
 $S(y,\psi) = S(y,\varphi)$.
\hfill{$\qed$}

\begin{lemma}\label{l2.1} For any $y \in \R^d$ and $\varphi\in\bN$, if
the set $S(y,\varphi)$ is bounded, then it satisfies 
\begin{align}\label{2.1}
\rho(y,\varphi^y)=\rho(y,(\varphi^y \cap S(y,\varphi))\cup \psi),
\quad \forall ~ \psi\subset\R^d\setminus S(y,\varphi),\psi\in\bN.
\end{align} 
\end{lemma}
{\em Proof.} 
Suppose $S(y,\varphi)$ is bounded and 
$\psi\in\bN$ with $ \psi\subset\R^d\setminus S(\varphi)$.
Since  $S(y,\varphi) = S(y,\varphi^y)$, we can and do
assume without loss of generality that $y \in \varphi$. 
Also, by Lemma \ref{stopsetlem}, it suffices to prove the result
in the case where
 $\varphi \subset S(y,\varphi)$, so we assume this too.

 Write 
$\varphi'$ for $\varphi \cup \psi$, and write $\rho(x)$ for
$\rho(x,\varphi)$ and $\rho'(x)$ for $\rho(x,\varphi')$.

Suppose $\rho(y) > \rho'(y)$. Let $x_1$ be a smaller grain-neighbour
of $y$ in $\varphi'$. Assume for now that $x_1 \in \varphi$.
Then by the hard-core property of the lilypond model
on $\varphi$, we have
$\rho(y)+\rho(x_1) \leq |x_1-y| = \rho'(y)+\rho'(x_1)$, and hence
 $\rho(x_1) < \rho'(x_1)$. 

Next let $x_2$ be a smaller grain-neighbour of $x_1$ in $\varphi$.
By the hard-core property of the lilypond model
on $\varphi'$, we have $\rho'(x_1) + \rho'(x_2) \leq |x_1-x_2|
= \rho(x_1)+ \rho(x_2)$, so 
$\rho'(x_2) < \rho(x_2)$. Let
$x_3$ be a smaller grain-neighbour of $x_2$ in $\varphi'$.
Assuming $x_3 \in \varphi $, using once more the hard-core property 
of the lilypond model on $\varphi$ 
yields $\rho(x_3) < \rho'(x_3)$. 

Continuing to alternate in this way, we get a sequence of
points $x_i$ satisfying
$$
\rho(y) > \rho'(y) \geq \rho'(x_1) > \rho(x_1) \geq \rho(x_2) > \rho'(x_2) \geq \rho'(x_3)
 > \rho(x_3) \geq   \cdots
$$ 
and terminating at $x_n$ if for some (odd) $n$ we
have $x_n \in \psi$.
We see from these inequalities that the possibly terminating
sequence  
$y, x_1, x_2, \ldots$
consists of distinct
points. Also, setting
$x_0 =y$ and $\rho_i =  \rho(x_i)$ and $\rho'_i =\rho'(x_i)$,
 we have $|x_i - x_{i-1}| = \rho'_i + \rho'_{i-1}$ 
for odd $i$ and 
$|x_i - x_{i-1}| = \rho_i + \rho_{i-1}$ 
for even $i$, and therefore from the above inequalities
$
|x_i - x_{i-1}| 
$ is nonincreasing (in fact, strictly decreasing) in $i$.

Thus the (possibly terminating) sequence $y,x_1,x_2,\ldots$
 forms a descending chain with $|x_1-y| \leq 
2 D(y,\varphi^y)$.
If the sequence terminates at some point $x_n=z$
with $z \in \psi$, then $z \in S(y,\varphi)$,
contradicting the assumption that $\psi \cap S(y,\varphi) = \emptyset$.
On the other hand, 
if the sequence $(x_i)$ does not terminate, then it 
 forms an infinite descending chain in $\varphi$,
 contradicting the assumption that  $S(y,\varphi)$ is bounded.

Thus if $\rho(y) > \rho'(y)$ we have derived a contradiction.
If $\rho(y) < \rho'(y)$ we argue similarly, this time
starting with $x_1$ a smaller grain-neighbour of $y$ in
$\varphi$. Again we end up with a contradiction. 
Thus we must have $\rho(y) = \rho'(y)$.
\hfill{$\qed$}

\vspace{0.3cm}
Given $\varphi\in\bN$, define 
\bea
R(\varphi) := \inf \{ r>0 :S(0,\varphi)\subset B_r(0) \},
\lbl{0611a}
\eea
with the convention $\inf(\emptyset) := +\infty$.

\begin{lemma} \lbl{LemRmeas}
The function $R: \bN \to [0,\infty]$ is Borel-measurable.
\end{lemma}
{\em Proof.} 
Let $t >0$. It suffices to prove $R^{-1}([0,t])$ is measurable,
i.e. in $\mathcal{N}$.
For $m \in \N$, let $\bN_m = \{ \varphi \in \bN: \varphi(B_t(0))=m \}$.
Then  by Lemma \ref{stopsetlem},
 $$
R^{-1}([0,t]) = \cup_{m =0}^\infty
\{ \varphi \in \bN_m : R(\varphi \cap B_t(0)) \leq t \}.  
$$
However, it is not hard to see that for all $m \in \N$,
the function 
$g: B_t(0)^m \to \{0,1\}$ given by
$$
g(x_1,\ldots,x_m) := \begin{cases} 
 {\bf 1}_{\{ R(\{x_1,\ldots,x_m\}) \leq t \} }
& {\rm if ~ } x_1,\ldots,x_m ~{\rm are~ distinct}
\\
0 & {\rm otherwise}
\end{cases}
$$
is measurable, and hence
$\{ \varphi \in \bN_m : R(\varphi \cap B_t(0)) \leq t \} \in \mathcal{N}$.
The result follows. $\qed$\\

By Lemma \ref{l2.1}, $R(\varphi)$ is a radius of stabilization
for $\rho(0,\varphi^0)$, i.e.
$\rho(0,\varphi^0)$ is unaffected by modifications to
$\varphi$ outside $B_{R(\varphi)}(0)$.
Also, $R(\Phi)$  is known to be almost
surely finite \cite{MR,DL05}.  The next result provides
another proof of this last fact, and more importantly
shows that the tail of the distribution
of $R(\Phi)$ decays  sub-exponentially.

\begin{lemma}\lbl{l2.2}
There is a constant $c_1 > 0$
such that for any $x \in \R^d$,
\begin{align}
\BP(R(\Phi)>t)\le c_1\exp (-c_1^{-1} t^{d/(d+1)}),\quad t>0. 
\label{Rtail}
\end{align}
\end{lemma}
{\em Proof:} 
Let $t\geq 1$.
 The first step is to bound the 
distance $D:=D(0,\Phi^0)$ from the origin to 
its nearest neighbour in $\Phi$ by a suitable power of $t$.
For the duration of this proof, define
 $$
\varepsilon := 1/(d+1) ; ~~~ K : = \max (2, (2 b_d e^2)^{1/(d+1)} ).
$$
Given $t \geq 1$,
set $u:= K^{-1}t^{\eps}$ and $\ell := \lfloor Kt^{1-\eps} \rfloor $.
Let $E$ be the event that there is a descending chain
$0,x_1,\ldots,x_\ell$ in $\Phi^{0}$ such that $K|x_1|\le t^\varepsilon$
(i.e., $|x_1| \leq u$).
Then we assert that
the event inclusion
\begin{align}\label{2.9}
\{ R(\Phi)>t \} \subset \{2KD>t^\varepsilon\} \cup E
\end{align}
holds.  To see this, assume on the contrary that $2KD\le t^\varepsilon$
and that $E$ does not occur. Then
any descending chain starting at $0$ and with its first link of
Euclidean length $\le 2D\le K^{-1}t^\varepsilon$ would have
at most $Kt^{1-\eps}-1$ links, so would end at a point of
Euclidean norm  at most  $u(Kt^{1-\varepsilon} -1) $. 
This would imply that $R(\Phi)\le uK t^{1-\eps} = t$.
Hence \eq{2.9} holds.

For any set $A$,
let $(A)_\ell$ denote the set of $\ell$-tuples 
 $(a_1,\ldots,a_\ell)\in A^\ell$ such that $a_1,\ldots,a_\ell$
are distinct. 
Then, using the $\ell$th order Palm-Mecke
 formula for the Poisson process (see e.g. Theorem 1.6 of \cite{Pe})
 similarly to Subsection 3.2 of \cite{DL05}, 
 we have
that
\bea
\BP(E)    & \leq & 
\BE \sum_{(x_1,\ldots,x_\ell) \in (\Phi)_\ell} {\bf 1}\{ u  \geq 
|x_1| \geq |x_2-x_1| \geq \ldots \geq |x_\ell-x_{\ell-1}| \}
\nonumber \\
& = &
\idotsint
\I\{u \ge |x_1| \ge |x_2-x_1|\ge \ldots\ge |x_{\ell}-x_{\ell-1}| \}dx_1\ldots 
dx_{\ell} 
\nonumber
\\
& = &
\idotsint
\I\{u \ge |x_1| \ge |y_2|\ge \ldots\ge |y_{\ell}| \}dx_1d y_2\ldots dy_{\ell} 
\nonumber \\
& =  &
\frac{ ( b_du^d )^{\ell} }{\ell!} .
\lbl{0312b}
\eea
Hence by Robbins' bound associated with 
Stirling's formula (see e.g. \cite{Feller1}),
 and the definition of $u$ and $\ell$, and the fact that $K \geq 2$,
\bean
\BP(E)
\leq  \frac{ b^\ell_d u^{\ell d} e^\ell  }{\ell^\ell (2 \pi)^{1/2} }
\leq \left( \frac{b_d t^{\eps d} e }{K^d(K t^{1-\eps} -1) } \right)^\ell
\leq \left( \frac{2 b_d t^{\eps d} e }{K^{d+1} t^{1-\eps}  } \right)^\ell.
\eean
By definition, $\eps d = 1 -\eps$, and $K^{d+1} \geq 2 b_d e^2$
so that 
\bean
\BP(E) \leq e^{-\ell} \leq
 \exp ( -(K/2) t^{d/(d+1)} ). 
\eean
Returning to \eqref{2.9} and using the definition of $\eps$ again,
we see that
\bea
\BP(R(\Phi) > t)\le
 \exp(-b_d(2K)^{-d}t^{d/(d+1)})+
 \exp (- (K/2) t^{d/(d+1)} ), 
\lbl{0610a}
\eea
and therefore \eq{Rtail} holds for suitably chosen $c_1$.
\hfill{$\qed$}\\

Next we extend
 Lemma \ref{l2.2} to the  family of
 binomial point processes 
\bea
\chi_{n,m} := \{n^{1/d} X_i: 1 \leq i \leq m\},
~~~~~ n,m \in \N.
\lbl{chidef}
\eea
 where 
$X_1,X_2,\ldots$ is a sequence of  independent random $d$-vectors
uniformly distributed over the convex set $W$ (which was 
introduced in Section \ref{secintro}).
\begin{lemma}\lbl{l2.2a} 
There is a constant $c_2>0$, dependent only
on the choice of  $W$, 
such that for all $n \in \N$ with $n \geq 4$,
and all $x \in  W_n$,
and all $m \in [n/2,3n/2]$,
\begin{align}\label{Rtail2}
\BP(R(-x +\X_{n,m-1})>t)\le c_2 \exp (-c_2^{-1} t^{d/(d+1)}),\quad t>0.
\end{align}
\end{lemma}
In proving this, and again later, we shall use
the fact that there is a constant $c_0>0$,
dependent only on $W$, such that  
\bea
c_0^{-1} r^d \leq V_d(W_n \cap B_{r}(x) ), ~~~ 
~ n \in \N, x \in W_n, 
  r \in (0, n^{1/d} \diam(W)]. 
\lbl{convlow}
\eea
To see this, take $u>0$ such that $B_{2u}(0) \subset W$.
For $x \in W $, and $0< s \leq u$,
the convex hull of $\{ x\} \cup B_{2u} (0)$ is contained
in $W$,  and the intersection of this with $B_{s}(x)$ contains 
the intersection of $B_{s}(x)$
 and a cone with apex at $ x$ and subtended angle
at least $\arcsin(u/\diam(W)) $,
 and hence has volume at least $c^{-1}s^d$ for some
constant $c$ independent of $x$ and $s$. Therefore
\bea
\inf \{ s^{-d} V_d(W \cap B_{s}(x) ): x \in W, s \in (0,u] \}
> 0.
\lbl{convlow0}
\eea
Since $s^{-d} V_d(W \cap B_{s}(x) ) $
is continuous in $(s,x)$, it is bounded away from zero
on $(s,x)$ in the compact set $ [u,\diam(W)] \times W$, and therefore
 \eq{convlow0} still holds
with the range of $s$ extended to $(0,\diam(W)]$, and then
\eq{convlow} follows by setting $x= n^{-1/d}y$ and scaling.  \\

{\em Proof of Lemma \ref{l2.2a}.}
By \eq{convlow},
for all $n \geq 4, x \in W_n,$ $m \geq n/2$  (so that in particular
$m- 1 \geq n/4$)
 and $r \in (0, n^{1/d} \diam(W)]$, we 
have
\bea
 \BP ( D(x,\X_{n,m-1} )
\geq r ) \leq (1- c_0^{-1}n^{-1}r^d)^{m-1}
 \leq (\exp( - c_0^{-1}n^{-1}r^d ))^{m-1}
\nonumber \\
\leq  \exp( - (4 c_0)^{-1}r^d),
\lbl{100104b}
\eea
and this  holds trivially for $r >  n^{1/d} \diam(W)$ as  well.

The proof of Lemma  \ref{l2.2a} now 
mostly follows that of Lemma \ref{l2.2}. There is a difference 
in the first term in the right side of  \eq{0610a},
where we now need to estimate the  probability that
the nearest point to $x$ in $\chi_{n,m-1}$
is at a distance greater than $(2K)^{-1}t^\eps$. For this we can use 
 \eq{100104b}.

Instead of the event $E$ featuring 
 in the proof of Lemma \ref{l2.2}, we now need to consider
 $E'$, defined to be the event that there is a descending chain
$(x,x_1,\ldots,x_{\ell})$ in $\chi_{n,m-1}^x$    with $|x_1 -x| \leq u$.
Instead of 
the estimate \eq{0312b},
setting $Y_i := n^{1/d} X_i$
 we now have
\bea
\BP(E') \leq 
\BE \sum_{(i_1 , i_2 , \ldots , i_{\ell } ) \in ( \{1,2,\ldots,m\})_\ell
} {\bf 1}\{ u  \geq 
| Y_{i_1}| \geq | Y_{i_2} - Y_{i_1}| \geq \cdots 
\geq | Y_{i_\ell}- Y_{i_{\ell-1}}| \}
\nonumber \\
= \frac{m!}{(m-\ell)!}
\int_{(n^{1/d} W)^{\ell}}
 {\bf 1}\{ u  \geq 
|x_1| \geq |x_2-x_1| \geq \cdots \geq |x_\ell-x_{\ell-1}| \}
n^{-\ell} dx_1 \ldots dx_\ell
\nonumber \\
 \le 
2^{\ell } 
\idotsint
\I\{u \ge |x_1| \ge |x_2-x_1|\ge \ldots\ge |x_{\ell}-x_{\ell-1}| \}dx_1\ldots 
dx_{\ell} 
\label{100115b}
\eea
because $m \leq 2n$.

With these changes,  we can complete the proof
by  following the proof of Lemma \ref{l2.2};
it is easy to modify the argument to allow
for the extra factor of $2^{\ell}$ in \eq{100115b} compared to
\eq{0312b}.
\hfill{$\qed$} 

\section{External stabilization}
\lbl{secstrong}
\eqco
In this section we introduce the notion
of a {\em fence}.  Loosely speaking, given an annulus in $\R^d$,
 a configuration $\varphi \in \N$
has a fence if there are enough points
in the annulus to guarantee
that no lilypond grain centred inside
the annulus can penetrate too far outside,
and no lilypond grain centred outside
the annulus can penetrate too far inside.
Combining this with with the notion
of a radius of stabilization $R(\varphi)$ as
already considered, we shall arrive at a stronger
{\em external} stabilization radius, denoted 
$\Rex(\varphi)$, with similar tail behaviour.
Loosely speaking, external stabilization means that
changes to $\varphi$ beyond distance 
 $\Rex(\varphi)$ do not affect the grains  near the origin,
and changes near the origin do not affect the 
grains beyond distance $\Rex(\varphi)$.

Let $x \in \R^d$ and $0<r <s$.
Let $B^o_r(x)$ be the open ball of radius $r$ centred at $x$, and
set $\partial B_s(x):= B_s(x) \setminus B^o_s(x)$, 
the boundary of $B_s(x)$.
Let $(e_1,e_2,e_3,\ldots)$ be an arbitrarily chosen
 sequence forming a countable dense
set in $\partial B_1(0)$. 
Let $k(s,r)$ 
be the smallest positive integer $k$ such that there exists
 an increasing sequence $(j_1,  j_2, \ldots, j_k)$ 
of positive integers, such that 
\bea
\partial B_s(0) \subset \cup_{i=1}^k
B^o_r(se_{j_i}  )
\lbl{forfence}
\eea 
 (such a $k$ exists by compactness).
Note that 
\bea
k(s,r) = k(s/r,1),
~~~ 
0 < r < s.
\lbl{0610c}
\eea
In other words, $k(s,r)$ depends on $(s,r)$ only through
the ratio of $s$ to $r$.

Let $x \in \R^d$ and $0< r <s$.  Setting $k=k(s,r)$,
let $j_1,\ldots,j_k$ be the first sequence of positive integers,
 according to the lexicographic ordering, such  that \eq{forfence}
holds. That is, for $0 \leq i \leq k-1$, having defined $j_1,\ldots,j_i$ let 
 $j_{i+1}$ be the first positive integer $j$
such that  $\{j_1\ldots,j_i,j\}$ can be extended to 
a set of cardinality $k$ in $\N$  which satisfies 
\eq{forfence}.

For any $x \in \R^d$,
if we set $y_i = x + s e_{j_i}$, then the  points $y_1,\ldots,y_k$
in $\partial B_s(x)$ satisfy
$$
\partial B_s(x) \subset \cup_{i=1}^k
B^o_r(y_i  ).
$$
With $y_1, \ldots, y_k$ defined thus, define 
\bea
F(x,s,r) := 
\cap_{i=1}^k \{ \varphi \in  \bN: \varphi(B^o_r(y_i) \setminus B_s(x)) 
\geq 2 \} .
\lbl{fencedef}
\eea
In the case $x=0$ we write simply $F(s,r)$ for $F(0,s,r)$.
We  think of $F(x,s,r)$ as a set of configurations $\varphi$ 
containing  a fence  in the annulus $B_{s+r}(x) \setminus B_s(x)$.
We formalize the fence property as follows:
\begin{lemma}
\lbl{trilem}
Suppose $0< r <s < \infty$ with $s \geq 2 r$, and let $x \in \R^d$.
Suppose $\varphi \in F(x,s,r)$.  Then 
for any $z \in \varphi \setminus B^o_s(x)$ we  have $B_{\rho(z,\varphi)}(z)
\cap B_{s-2r}(x) = \emptyset$, and 
for any $y \in \varphi \cap  B_s(x)$ we  have $B_{\rho(y,\varphi)}(y)
\subset B^o_{s+2r}(x) $.
\end{lemma}
{\em Proof.}
Let $z \in \varphi \setminus B^o_{s}(x)$.
The line segment from $z$ to $x$ includes a point $w$
 in $\partial B_{s}(x)$.
Then $w$ lies in at least one of the 
balls $B^o_r(y_i)$, 
and since we assume $\varphi \in F(x,s,r)$,
there exists $u \in \varphi \cap B_r^o(y_i) \setminus \{z\}$ such
that $|u-w| <2r$. 
Hence by the triangle inequality,
$|z-u| < |z-w| + 2r$. Therefore,
 since $\rho(z,\varphi) \leq |z-u|$, 
the grain $B_{\rho(z,\varphi)}(z)$ does not intersect $B_{s-2r}(x)$,
as asserted.

Let $y \in \varphi \cap  B_{s}(x)$.
Take $v$ in $\partial B_{s}(x)$,
such that $y$ lies in 
the line segment from $v$ to $x$.
Since 
 $\varphi \in F(x,s,r)$,
there exists $t \in \varphi  \setminus \{y\}$ such
that $|t-v| <2r$. 
Hence by the triangle inequality,
$\rho(y, \varphi) \leq  |y-t| < |y-v| + 2r$, so that
 $B_{\rho(y,\varphi)}(y) \subset B^o_{s+2r}(x)$,
as asserted.
\hfill{$\qed$} \\

Recall that $W\subset\R^d$ is  convex and compact with $V_d(W)=1$,
 containing an open neighborhood of the origin, and
 $W_n : = n^{1/d} W_n$ for $n \in \N$.
For proving results on point processes in
$W_n$, 
we introduce some further notation.
Given $0< r < s$, given $n \in \N$ and $x \in W_n$, we shall define
$F_n(x,s,r)$ similarly to $F(x,s,r)$ but now with the fence 
 involving only regions intersecting  $W_n$. 
First we define the points $y_i = x+ se_{j_i}$ for
$ 1 \le i \leq k (s,r/2)$, similarly to the points $y_i$  in
the definition of $F(x,s,r)$, such that
\bea
\partial  B_s(x) \subset \cup_{i=1}^{k(s,r/2)} B^o_{r/2}(y_i).
\lbl{0604a}
\eea
Let $z_1,\ldots,z_{k'}$ be those $y_i$
such that $B_{r/2}(y_i) \cap W_n \neq \emptyset$. Note that
$k' \leq k(s,r/2)$. Set
\bea
F_n (x,s,r) := \cap_{i=1}^{k'} \{ \varphi \in \bN : \varphi( B_r^o(z_i)) \geq 2 
\} \cap \{\varphi \in \bN: \varphi \subset W_n\} .
\lbl{100115a}
\eea 

\begin{lemma}
\lbl{tri2lem}
Suppose $\varphi \in \bN$ and $\psi \in \bN$ (possibly
with $\psi = \varphi$). Let $r,s,t >0$ with $s \geq 2r$ and
$t \geq s+4r$, and $x \in \R^d$.
  Suppose $y \in \varphi \cap B_s(x)$ and $z \in \psi \setminus B_t(x)$ 
Then:

(i) If   $\varphi \in F(x,s,r)$ and $\psi \in F(x,t,r)$
 then
$|y-z| > \rho(y,\varphi) + \rho(z,\psi)$.

(ii) If   $\varphi \in F_n(x,s,r)$ and $\psi \in F_n(x,t,r)$
 then
  $|y-z| > \rho(y,\varphi) + \rho(z,\psi)$.
\end{lemma}
{\em Proof.}
To prove (i), suppose
   $\varphi \in F(x,s,r)$ and $\psi \in F(x,t,r)$.
 Then by  Lemma \ref{trilem}, since $t-s \geq 4r$ we have 
$B_{\rho(y,\varphi)} \subset B_{s+2r}^o(x)$ and 
$B_{\rho(z,\psi)} \cap B_{t-2r}(x) = \emptyset$,
so that
$B_{\rho(y,\varphi)} (y) \cap B_{\rho(z,\psi)} (z)  = \emptyset$,
establishing part (i).

To prove (ii), suppose instead that   $\varphi \in F_n(x,s,r)$ and
 $\psi \in F_n(x,t,r)$.
Let $u$ and $v$ be the points on the line segment
$yz$ such that $|u -x| = s$ and $|v-x|= t$.
Since $W_n$ is convex, $u$ and $v$ are in $W_n$.
Since  $u \in \partial B_s(x)$,
by \eq{0604a} in the definition of $F_n(x,s,r)$ we have $u \in B_{r/2}(y_i)$
for some $i$, but then since  also $u \in W_n$ we have
$y_i = z_j$ for some $j \leq k'$. Hence since
$\varphi \in F_n(x,s,r)$, by \eq{100115a},
there exists $w \in \varphi \cap B_{r}^o(z_j) \setminus \{y\}$,
and   by the triangle inequality 
$$
|w-y | \leq |u-y| + |z_j-u| + |w-z_j| \leq |u-y| + 3r/2.
$$ 
Similarly, since $v \in W_n \cap \partial B_t(x)$ and
$\psi \in F_n(x,t,r)$ we can find $w' \in \psi \setminus \{z\}$ such that
$
|w'-z | \leq |v-z|  + 3r/2.
$ 
Hence by the hard-core property,
$$
\rho(y,\varphi) + \rho(z,\psi) \leq |w-y| + |w'-z| \leq
|u-y| + |v-z| + 3r,
$$
whereas since $y,u,v,z$ are collinear and
 $|u-v| \geq t-s \geq 4r$, we have
$$
|y-z|  \geq |u-y| + |v-z| + 4r,
$$
and part (ii) follows. \hfill{$\qed$} \\

We now give some probability estimates for the
point sets $F(x,s,r)$ and $F_n(x,s,r)$.

\begin{lemma}
\lbl{lemFprob}
 There exists $c_3 \in(0,\infty)$ such that
if $0<r<2r <s$ and $x \in \R^d$, then:

(i) $\BP (\Phi \notin F(x,s,r)) \leq k(s,r) c_3 \exp(- c_3^{-1} r^d)$.

(ii) If $m,n \in \N$ and $n \geq 4$ and $m \in [n/2,3n/2]$ and
$x \in W_n$, then 
 $\BP (\chi_{n,m}  \notin F_n(x,s,r)) \leq k(s,r) c_3 \exp(- c_3^{-1} r^d)$.
\end{lemma}

{\em Proof.}
(i)
With $y_1, \ldots, y_k$ as in the definition \eq{fencedef}
of $F(x,s,r)$,
note that $V_d(B_r^o(y_i) \setminus B_s(x)) \geq (b_d/2) r^d$,
so by subadditivity of measure,
\bea
\BP (\Phi \notin F(x,s,r)) \leq
k(s,r)
 (1+ (b_d/2) r^d) \exp( - (b_d/2) r^d ),
\lbl{0311d}
\eea 
and part (i) follows.

For part (ii), first we  claim that there is a constant
$c >0$, independent of $n$ and $x$, such that
if $n \geq 4$ and $0 \leq r \leq n^{1/d} \diam(W)$, then
with $k'$ and $z_i$ as defined just before  \eq{100115a} we have
 \bea
V_d(B_r(z_i) \cap W_n) \geq c^{-1} r^d, ~~~
1 \leq i \leq k'.
\lbl{0116a}
\eea
Indeed, given $i\leq k'$,
 if we choose $y \in B_{r/2}(z_i) \cap W_n$, 
then $B_{r/2}(y) \subset B_r(z_i)$
so the claim follows from \eq{convlow}.

By \eq{0116a}, subadditivity and the fact that the binomial
distribution is stochastically increasing in the
success probability, we have
\begin{align*}
\BP (\chi_{n,m} \notin F_n(x,s,r)) \le
k(s,r) [(1-r^d/(nc))^m+m(r^d/(nc)) (1- r^d/(nc))^{m-1}].
\end{align*} 
Hence by 
 the inequality $1-t \leq e^{-t}$, 
there is a further constant $c'$ such that
for all $n \geq 4$ and $m \in [n/2,3n/2]$,
 (so that in particular $m-1 \geq n/4$ and $m \leq 2n$),
 for all $x \in W_n$ and $0<r< s$ with  $ r \leq  n^{1/d} \diam(W)$,
\bea
\BP [\X_{n,m} \notin F_n(x,s,r)] \leq
k(s,r) (1+ m   r^d/(cn)) (1-  r^d/(c n))^{m-1} \nonumber \\
\leq k (s,r) (1+(2/c)r^d) \exp (- (4c)^{-1} r^d) \nonumber
 \\
\leq k(s,r) c' \exp(-(c')^{-1} r^d).
\lbl{0311d2}
\eea 
Moreover, if $r > n^{1/d} \diam (W)$ then $s-r > n^{1/d} \diam (W)$
so $k'=0$, and then
 trivially \eq{0311d2} still holds. 
This gives us part (ii). \hfill{$\qed$} \\

For $\varphi \in \bN$ and $r >0$, define the set
\bea
S_r^*(\varphi) := \cup_{x \in \varphi \cap B_{7r}(0) \setminus B_{2r}(0)}
S(x,\varphi).
\lbl{Sdef}
\eea

\begin{lemma}\lbl{fencelem} Let $r >0$.  Suppose
$\varphi \in F(2r,r/2) \cap F(4r, r/2) \cap F(7r,r/2)$, and
suppose $S_r^*(\varphi) $ is bounded. 
Let $\varphi_{{\rm in}} \in \bN , \varphi_{{\rm out}} \in \bN$ be
 such that 
$$
\varphi_{{\rm in}} \subset B_{2r}(0) \setminus S_r^*(\varphi); ~~~~
\varphi_{{\rm out}} \subset  \R^d \setminus (B_{7r}(0) \cup  S_r^*(\varphi)).
$$
Then
\bea
\rho (x, \varphi \cup \varphi_{{\rm in}} ) = \rho(x,\varphi 
\cup \varphi_{{\rm in}}
\cup \varphi_{{\rm out}}), ~~~ 
~~~
 x \in
 \varphi_{{\rm in}}
 \cup
 (\varphi \cap B_{7r}(0))
; 
\lbl{rexin}
\\
\rho (x, \varphi \cup \varphi_{{\rm out}} ) = \rho(x,\varphi 
\cup \varphi_{{\rm in}}
\cup \varphi_{{\rm out}}), 
~~~ 
~~~ 
x \in
 \varphi_{{\rm out}}
 \cup
 (\varphi \setminus B_{2r}(0))
. 
\lbl{rexout}
\eea
\end{lemma}
{\em Proof.}
First note that by  Lemma \ref{l2.1},
\bea
\rho(x, \varphi \cup \varphi_{{\rm in}}) =
\rho(x, \varphi \cup \varphi_{{\rm out}}) =
\rho(x, \varphi ), ~~~  ~~~
 x \in \varphi \cap B_{7r}(0) \setminus 
B_{2r}(0). 
\eea
Hence, 
we can and do consistently define $\rho'(x)$
for all $x \in \varphi \cup \varphi_{{\rm in}} \cup \varphi_{{\rm out}} $,
 by
\bea
\rho'(x) = \rho(x,\varphi \cup \varphi_{{\rm in}}), ~~~ x 
\in \varphi_{{\rm in}} \cup ( \varphi \cap B_{7r}(0));
\lbl{rexin2}
\\
\rho'(x) = \rho(x,\varphi \cup \varphi_{{\rm out}}), ~~~ 
x \in \varphi_{{\rm out}} \cup (\varphi \setminus B_{2r}(0)) .
\lbl{rexout2}
\eea
Assign grain radius $\rho'(x)$ to each 
 $x \in \varphi \cup \varphi_{{\rm in}} \cup \varphi_{{\rm out}} $.
We assert that with  grain radii assigned in this way,
 each  
$x \in \varphi \cup \varphi_{{\rm in}} \cup \varphi_{{\rm out}}$
 has a smaller grain-neighbour
 in $ \varphi \cup \varphi_{{\rm in}} \cup \varphi_{{\rm out}} $,

To verify this assertion, first suppose 
 $x \in B_{4r}(0)$. Then by the defining properties
of the lilypond model, there exists $y \in  \varphi \cup 
\varphi_{{\rm in}}$ such that $y$ is a smaller grain-neighbour
of $x$ under the radii $\rho(\cdot, \varphi \cup \varphi_{\rm in})$.
Moreover,
since $\varphi \in F(4r,r/2) \cap F(7r,r/2)$, 
if  $y \notin B_{7r}(0)$ then 
by Lemma \ref{tri2lem},
we would have
$|x-y| > \rho(x,\varphi \cup \varphi_{{\rm in}} )
+ \rho(y,\varphi \cup \varphi_{{\rm in}})$,
contradicting the statement that $y$ is a grain-neighbour of $x$.
Therefore $y \in B_{7r}(0)$, so by \eq{rexin2}
we have $\rho'(y) = \rho(y, \varphi \cup \varphi_{\rm in})$ 
(and likewise for $x$). Therefore $y$ is also a smaller
grain-neighbour of $x$ using the radii $\rho'(\cdot)$ as asserted.

Now suppose instead that 
 $x \notin B_{4r}(0)$. Then there exists $y \in  \varphi \cup 
\varphi_{{\rm out}}$ such that $y$ is a smaller grain-neighbour
of $x$ under the radii $\rho(\cdot, \varphi \cup \varphi_{\rm out})$.
Moreover,
since $\varphi \in F(4r,r/2) \cap F(2r,r/2) $,
if $y \in B_{2r}(0)$ then 
 by Lemma \ref{tri2lem},
we would have
$|x-y| > \rho(x, \varphi \cup \varphi_{{\rm out}}) +
 \rho(y, \varphi \cup \varphi_{{\rm out}}) $
contradicting the statement that $y$ is a grain-neighbour of $x$.
Therefore $y \notin B_{2r}(0)$, so by \eq{rexout2}
we have $\rho'(y) = \rho(y, \varphi \cup \varphi_{\rm out})$ 
(and likewise for $x$). Therefore $y$ is also a smaller
grain-neighbour of $x$ using the radii $\rho'(\cdot)$ as asserted.

We shall show that the radii  $\rho'(x), x \in
 \varphi \cup \varphi_{{\rm in}} \cup \varphi_{{\rm out}},$
 have the hard-core property.
This will suffice to  give
 $\rho'(x) = \rho(x, \varphi \cup \varphi_{{\rm in}} \cup \varphi_{{\rm out}})$
for all $x \in  \varphi \cup \varphi_{{\rm in}} \cup \varphi_{{\rm out}}$ 
as required, because
as already mentioned, 
 the lilypond model
is the {\em unique} set of radii satisfying the hard-core
and smaller grain-neighbour properties \cite{HL06}.

Let $x,y$ be 
 distinct  elements of 
$ \varphi \cup \varphi_{{\rm in}} \cup \varphi_{{\rm out}}$, with
$|x| \leq |y|$. We consider separately the case with $|y| \leq 7r$,
the case with $|x| > 2r$, and the case with $|x| \leq 2r $ and $|y| > 7r$. These
three cases cover all possibilities.

In the first case with $|y| \leq 7r$, both $x$ and $y$ are   in
$\varphi_{{\rm in}} \cup (\varphi \cap B_{7r}(0))$. By \eq{rexin2}
and the hard-core property of the lilypond model
 on $\varphi \cup \varphi_{{\rm in}}$, we have
$\rho'(x) + \rho'(y) \leq |x-y|$.  

In the second case with $|x| > 2r$, both $x$ and $y$ are   in
$\varphi_{{\rm out}} \cup (\varphi \setminus B_{2r}(0))$. By \eq{rexout2}
and the hard-core property of the lilypond model on
 $\varphi \cup \varphi_{{\rm out}}$, 
we have $\rho'(x) + \rho'(y) \leq |x-y|$.  

Now consider the third case with
 $|x| \leq 2r $ and $|y| > 7r$. In this case we have
$x \in \varphi_{{\rm in}} \cup (\varphi \cap B_{2r}(0))$
 and $y \in \varphi_{{\rm out}} \cup (\varphi \setminus B_{7r}(0))$, so
by \eq{rexin2} and \eq{rexout2},
 the assumption that $\varphi \in F(2r,r/2) \cap F(7r,r/2)$,
and Lemma \ref{tri2lem}, we have that
$\rho'(x) + \rho'(y) < |x-y|$.

Hence the radii  $\rho'(x), x \in
 \varphi \cup \varphi_{{\rm in}} \cup \varphi_{{\rm out}}$
 have the hard-core property as required.
\hfill{$\qed$} \\

We have a similar result to Lemma \ref{fencelem} 
in the case of point processes in $W_n$.

\begin{lemma}\lbl{fencenlem} Let $r >0$ and $z \in \R^d$.
  Suppose 
$\varphi \in F_n(z,2r,r/2) \cap F_n(z,4r, r/2) \cap F_n(z,7r,r/2)$,
with $\varphi(\R^d) \geq 2$.
Let $\varphi_{{\rm in}} \in \bN , \varphi_{{\rm out}} \in \bN$ be
 such that 
$$
\varphi_{{\rm in}} \subset W_n \cap B_{2r}(x) \setminus S_r^*(\varphi); ~~~~
\varphi_{{\rm out}} \subset  W_n \setminus (B_{7r}(x) \cup  S_r^*(\varphi)).
$$
Then
\bean
\rho (x, \varphi \cup \varphi_{{\rm in}} ) = \rho(x,\varphi 
\cup \varphi_{{\rm in}}
\cup \varphi_{{\rm out}}), ~~~ 
~~~
 x \in
 \varphi_{{\rm in}}
 \cup
 (\varphi \cap B_{7r}(z))
; 
\\
\rho (x, \varphi \cup \varphi_{{\rm out}} ) = \rho(x,\varphi 
\cup \varphi_{{\rm in}}
\cup \varphi_{{\rm out}}), 
~~~ 
~~~ 
x \in
 \varphi_{{\rm out}}
 \cup
 (\varphi \setminus B_{2r}(z))
. 
\eean
\end{lemma}
{\em Proof.} The proof is just the same
as for Lemma \ref{fencelem}, only using part (ii)
instead of part (i) 
of Lemma \ref{tri2lem}. \hfill{$\qed$} \\ 

We now define $\bU$ to be the set of all $\varphi \in \bN$ such that  
every point of $\varphi$
has a {\em unique smaller grain neighbour} under the lilypond
model based on $\varphi$, or $\varphi(\R^d) \leq 2$.

Recall the definition \eq{0611a} of
 $R(\varphi)$, $\varphi \in \bN$.
  For $x \in \R^d$, and $r >0$, $n \in \N$ define subsets $E_r(x), U_r(x),
G_r(x)$ and
 $ G_{n,r}(x)$ of $\bN$ by
 \bean
E_r(x)  & : = &
\{ \varphi \in \bN: R(-y+ \varphi) < r, ~
 y \in \varphi \cap B_{8r}(x)
 \}; 
 \\
U_r(x)  &:=  &\{\varphi \in \bN: \varphi \cap B_{9r}(x) \in \bU \}, 
\eean
and
\bea
G_r(x)  &: = & E_r(x) \cap U_r(x) \cap \cap_{j=1}^8 F(x,jr,r/2) ;
\lbl{Grdef}
\\
G_{n,r}(x) & : = & 
 E_r(x) \cap U_r(x) \cap \cap_{j=1}^8 F_n(x,jr,r/2).
\lbl{Gnrdef}
\eea
%


For $\varphi \in \bN$, we now define our radius of external
 stabilization $\Rex(\varphi)$ by
\bea
\Rex(\varphi) := 9 \min\{r \in \N : \varphi \in G_r(0)
\},
\label{0605a}
\eea
with $\min(\emptyset)$ taken to be $+\infty$. 
The next result, in which  we write $B_r$ for $B_r(0)$, 
shows that $\Rex$ has the
  external stabilization property.
\begin{lemma}
\lbl{Rexlem}
Suppose $r>0$.
Then (i)
if  $\varphi \in G_r$, then
 for any $\psi \in \bN$ with $\psi (B_{8r})=0$,
we have:
\bea
  \rho( y, (\varphi \cap B_{8r}) \cup \psi) 
=  \rho( y, (\varphi^0 \cap B_{8r}) \cup \psi), 
~~ 
~~~ y \in \psi \cup (\varphi \cap B_{8r} \setminus B_{2r}); 
\lbl{Rex3}
\\
 \rho( x, (\varphi \cap B_{8r}) \cup \psi) 
=
 \rho( x, \varphi \cap B_{8r} ),  
~~ 
~~ ~ x \in \varphi \cap B_{7r} ;
\lbl{Rex1}
\\
 \rho( x, (\varphi^0 \cap B_{8r}) \cup \psi) 
=
 \rho( x, \varphi^0 \cap B_{8r} ), 
~~ 
~~ 
 ~ x \in \varphi^0 \cap B_{7r} .
\lbl{Rex2}
\eea
Also, (ii) if $n \in \N$, and $x \in W_n$, and $\varphi \in G_{n,r}(x)$,
then
\bea
  \rho( y, \varphi) =  \rho( y, \varphi^x ), 
~~ ~~ 
~ y \in \varphi \setminus B_{2r}(x) 
\lbl{binRex3}
\eea
\end{lemma}
{\em Proof.}
(i) Suppose $\varphi \in G_r$. Then
$\varphi \in F(2r,r/2) \cap F(4r,r/2) \cap F(7r,r/2)$, and moreover 
$\varphi \in E_{r}(0)$ so
the set $S_r^*(\varphi)$ defined by
 \eq{Sdef} is contained in $B_{8r} \setminus B_r$.
Also $S_r^*(\varphi)= S_r^*(\varphi \cap B_{8r}(0))$ by Lemma
\ref{stopsetlem}.
Therefore we can apply Lemma \ref{fencelem} to $\varphi \cap B_{8r}(0)$; taking
$\varphi_{{\rm out}} = \psi$ and $\varphi_{{\rm in}} = \{0\}$,
we obtain \eq{Rex3} from \eq{rexout} and
\eq{Rex2} from \eq{rexin}, and
 taking 
$\varphi_{{\rm in}} = \emptyset$ we obtain 
 \eq{Rex1} from \eq{rexin},
 completing the proof of part (i).

Part (ii) is proved by a similar argument,
using Lemma \ref{fencenlem} instead of
Lemma \ref{fencelem}.
\hfill{$\qed$} \\

The next result gives us tail bounds on $\Rex(\Phi)$.

\begin{lemma} \lbl{Rexfinlem} There is a constant $c_4 \in (0,\infty)$ such that
for all $r >0$,
\bea
\lbl{0611e}
\BP ( \Phi \notin G_r(0) ) \leq c_4 \exp(-c_4^{-1} r^{d/(d+1)} )
\eea
and for $n  \geq 4$ and $m \in [n/2,3n/2]$ and $y \in W_n$,
\bea
\lbl{0611d}
\BP (\chi_{n,m} \notin G_{n,r}(y)  ) \leq c_4 \exp(-c_4^{-1} r^{d/(d+1)} ).
\eea
\end{lemma}
{\em Proof.}
By Lemma 5.1 of \cite{DL05}, for all $r  >0$
we have $\BP (\Phi \in U_r(0)) =1$. 
Also, by the Palm-Mecke equation, see e.g. Theorem 9.22 of
\cite{Pe}, and \eq{Rtail},
\bea
\BP (\Phi \notin E_{r}(0))
\leq \BE \sum_{x \in \Phi \cap B_{8r}(0)} {\bf 1}\{
R(-x +  \Phi) \geq r  \} 
\nonumber \\ = 
b_d (8r)^d \BP (R(\Phi) \geq r)
\leq b_d (8r)^d c_1 \exp (-c_1^{-1}r^{d/(d+1)}).
\lbl{0305b}
\eea
Also, by Lemma \ref{lemFprob} (i)  and \eq{0610c},
\bean
 \BP ( \Phi \notin \cap_{j=1}^8 F(jr,r/2) )
\leq  \sum_{j=1}^8 k(2j,1) c_3 \exp(- c_3^{-1} r^d) 
\eean
and combined with \eq{0305b} this gives \eq{0611e} for a suitable choice
of $c_4$.

Next, suppose $n \geq 4$ and $n/2 \leq m \leq 3n/2$, and
$y \in W_n$. Then by Lemma \ref{l2.2a},
\bea
\BP ( \chi_{n,m} \notin E_r(y) ) \leq \frac{m}{n} \int_{B_{6r}(y)}
\BP (R(-x + \chi_{n,m-1} ) > r) dx 
\nonumber \\
\leq (3/2) c_2 \exp (- c_2^{-1} r^{d/(d+1)} ). 
\lbl{0610d}
\eea
Combined with Lemma  \ref{lemFprob} (ii) 
 and \eq{0610c}, this gives us \eq{0611d}.
\hfill{$\qed$}

\section{Sub-exponential decay}
\lbl{secdecay}
\eqco

In this section we give tail bounds on the distribution
of the size of the component of the Poisson lilypond model 
containing a typical Poisson point, with `size'
measured either by cardinality, or by metric
diameter, or by volume.

For $\varphi \in \bN$, define as follows the  directed graph
$\G(\varphi)=(\varphi,E(\varphi))$ with vertex set $\varphi$
and edge set $E(\varphi)$. A pair $(x,y)$
is in $E(\varphi)$ if $y$ is a smaller grain-neighbour
of $x$. 
Let $\G^*(\varphi)$ denote the associated undirected graph.
For $x \in \varphi$,   let $C(x,\varphi)$ denote the cluster at $x$, that 
is, the set of points of $\varphi$ that are connected
to $x$  by a path in the undirected graph $\G^*(\varphi)$.
Let $C'(x,\varphi)$ denote the union of lilypond grains centred
at points of $C(x,\varphi)$, i.e., the connected
component containing $x$ of the
 set $\cup_{x \in \varphi} B_{\rho(x,\varphi)}(x)$.

\begin{theorem}
\lbl{thtail}
There are strictly positive constants $c_5$, $c_6,$ $c_7$
 such that
\bea
\label{0311e}
c_5^{-1} \exp(- c_5 r^d)  \leq & \BP (  \diam(C'(0,\Phi^0)) \geq r   ) & 
\leq c_5 \exp( - c_5^{-1} r^{d/(d+1)} ), ~~ r >0  ;
\\
\label{voltail}
c_6^{-1} \exp(- c_6 t)  \leq & \BP (  V_d (C'(0,\Phi^0)) \geq t   ) 
& \leq c_6 \exp( - c_6^{-1} t^{1/(d+1)} ), ~~ t > 0  ;
\\
\label{c1}
c_7^{-1} \exp(- c_7 n^{2})  \le & \BP(\card C(0,\Phi^0)\ge n)
& \le c_7 \exp (-c_7^{-1} n^{d/(d+1)} ) ,\quad n > 0.
\eea
\end{theorem}
Note that  in each of \eq{0311e}, \eq{voltail} and \eq{c1}
the power of $r$, $t$ or $n$ in the exponent is different in the lower bound
than in the upper bound. It is an open problem
to make these bounds sharper. 

Sharper  bounds are available in the analogous setting for 
lattice and continuum percolation. Consider for example
the geometric graph on $\Phi^0$, with each pair of points
connected by an edge if and only if
the distance between them is less than a constant $r^*$, with
$r^*$ chosen to be subcritical. Then results
like \eq{c1} and \eq{0311e} hold with exponents of the form $c\cdot n$,
respectively $c \cdot r$, 
in both the upper and lower bound, though not necessarily both
with the same $c$ (see Section 10.1 of 
\cite{Pe}). 
 A result like (\ref{0311e}) holds for 
subcritical lattice percolation
with exponents of the form $c\cdot r$;
see (6.10) of \cite{Grimm}.
This implies a bound like \eq{voltail} for the
Boolean model associated with the subcritical geometric graph
just mentioned,  with $c \cdot t$ in the exponent.
A similar lower bound also holds.  \\

The next lemma will be used in proving Theorem \ref{thtail},
and again later. 
Recall the definitions \eq{Grdef} and \eq{Gnrdef}
of $G_r(x)$ and $G_{n,r}(x)$ respectively.

\begin{lemma}
\lbl{prefroglem}
Let $r >0$ and $x \in \R^d$.

 (i) If $ \varphi \in G_r(x)  $, then 
$C(y,\varphi) \subset B_{5r}(x)$ for all
$y \in \varphi \cap B_{3r}(x)$, and 
\bea
\bigcup_{y \in \varphi \cap B_{3r}(x)} C'(y,\varphi) \subset B_{6r}(x).
\label{0611c}
\eea

(ii) If $n \in \bN$, and $x \in W_n$, and $\varphi \in G_{n,r}(x)$, then 
$C(y,\varphi) \subset B_{5r}(x)$. 
\end{lemma}
{\em Proof.} 
(i) First, we assert  that  each $y \in \varphi \cap B_{5r}(x)$ has
a unique smaller grain neighbour in $\varphi$.
Indeed, $y$ does not have any grain-neighbour
(in either  $\varphi$ or $\varphi \cap B_{9r}(x)$)
outside $B_{7r}(x)$, by Lemma \ref{tri2lem}
 because $\varphi \in F(x,5r,r/2) \cap 
F(x,7r,r/2)$. Also, 
 $\rho(u,\varphi) = \rho(u,\varphi \cap B_{9r}(9x))$
for all $u \in \varphi \cap B_{8r}(x)$ by Lemmas 
 \ref{stopsetlem} and \ref{l2.1},
because $\varphi \in E_r(x)$. Hence,
the unique smaller grain neighbour
of $y$ in $\varphi \cap B_{9r}(x)$
(which it has because $\varphi \in G_r(x) \subset U_r(x)$) is also
 its unique smaller grain neighbour
 in $\varphi $, justifying the assertion.
 
Hence, in the graph $\G(\varphi)$, each vertex
inside $B_{5r}(x)$
has an out-degree of 1. 
This implies that for any
path in $\G^{*}(\varphi)$ starting inside $B_{3r}(x)$
and ending outside $B_{5r}(x)$ but with all vertices
except the last inside $B_{5r}(x)$, if the direction of the edges
in $\G(\varphi)$ is taken into consideration the
path can reverse its direction at most once.
That is, such a path must consist of a
 directed path in the forward direction (possibly of zero length),
 followed by a  directed path in the reverse direction (also 
possibly of zero length). 
Hence, if $y \in \varphi \cap B_{3r}(x)$ and
$C(y,\varphi) $ is not contained in $B_{5r}(x)$, then
either there is a descending grain-chain in $\varphi$ starting 
inside $B_{3r}(x)$
 and ending outside $B_{4r}(x)$, or there is a descending
grain-chain  in $\varphi$ starting outside
$B_{5r}(x)$ and ending inside $B_{4r}(x)$. 

Since $\varphi \in E_r(x)$, we have no descending chain
starting inside $B_{3r}(x)$ and ending outside $B_{4r}(x)$ 
or starting in $B_{6r}(x)\setminus B_{5r}(x)$ and ending inside $B_{4r}(x)$.
Moreover, since $\varphi \in F(x,4r,r/2) \cap F(x,6r,r/2)$,
by Lemma \ref{tri2lem}
there is no  edge of $\G(\varphi)$
with
one endpoint
 outside $B_{6r}(x)$
and the other
endpoint
 inside  $B_{4r}(x)$.
This shows that $C(y,\varphi) \subset B_{5r}(x)$,
and since also $\varphi \in F(x,5r,r/2)$, by Lemma
\ref{trilem}  we have   \eq{0611c}.

The proof of (ii) is similar. 
\hfill{$\qed$} \\

{\em Proof of  \eq{0311e}.} 
By Lemmas 
 \ref{prefroglem}
 and
\ref{Rexfinlem}, we have for
$r >0$ that
\bea
 \BP (\{C'(0,\Phi^0) \subset B_{6r}
(0)\}^c) \leq  \BP(\Phi \notin G_r(0)) \leq c_4 \exp(-c_4^{-1} r^{d/(d+1)}),
\lbl{0612a}
\eea  
and the upper bound in \eq{0311e}   follows.

For the lower bound, observe that for all $r >0$,
we have 
\bean
\BP [ \diam (C'(0,\Phi^0) ) \geq r ] \geq
\BP [\Phi(B_{2r}(0))=0]
=
 \exp (- b_d (2r)^d).
\hfill{\qed}
\eean
\vspace{0.5em}

{\em Proof of \eq{c1}.}
 Choosing $r $ so that $2b_d (6r)^d = n$, 
by \eq{0612a} and a standard Chernoff-type 
tail estimate for the Poisson
distribution (see e.g. Lemma 1.2 of \cite{Pe})
there are constants $c,c'$ such that
 for  $n \geq 1$ 
we have
\bean
\BP [ \card ( C(0,\Phi^0)) \ge n ] \leq 
 c \exp(-c^{-1} r^{d/(d+1)} ) 
+ \BP [ \Phi^0 ( B_{6r}(0)) \geq 2 b_d (6r)^d ]
\\
\leq c' \exp (- (1/c') n^{1/(d+1)}) 
\eean 
and the upper bound in \eq{c1} follows.

Let $e$ be the unit vector $(1,0,\ldots,0)$ in $\R^d$.
For $i \in \N$,  let $B_i$ be the
 closed ball in $\R^d$,  centred on $9^{-i}e$ and
having radius $9^{-i-2}$. Since  $9^{-(i-1)} - 9^{-i} =
8(9^{-i})$,  for all $i \geq 2$, $x \in B_i$ and $y \in B_{i-1}$ we have
$7 (9^{-i}) \leq |x-y| \leq 9^{-(i-1)}$ and also
 $|x| \leq 10(9^{-i-1})$.

Observe that there is a positive constant $c$ such that
for all  $n \in\bN$ we have 
\bea
\BP ( \{ \Phi ( B_{9}(0)) =n \} \cap \cap_{i=1}^n \{ \Phi(B_i) =1\})
 = \left( \prod_{i=1}^n (b_d 9^{-(i+2)d}) \right) \exp(-b_d 9^d) 
\nonumber \\
\geq c^{-1} \exp(-c n^2). 
\label{LBprob}
\eea

If the event inside  the left hand side of \eq{LBprob} occurs,
 then labelling the point of $\Phi \cap B_i$
as $x_i$,  each point $x_i$ for $1 \leq i \leq n$
 has a smaller grain-neighbour in $\Phi^0$ to its left in
the collection $\{0,x_1,x_2,\ldots,x_n\}$. Indeed,
every point in this collection to the right 
of $x_i$, and also every point of $\Phi \setminus B_9(0)$,
is distant more than $2 D(x_i,\Phi^0)$
from $x_i$, so cannot be its  smaller grain neighbour.

Therefore, if the event inside  the left hand side of \eq{LBprob} occurs,
then $x_1,\ldots,x_n$ are all in $C(0,\Phi^0)$,
so that in this case $\card C(0,\Phi^0) \geq n+1$.
Hence, the lower bound in \eq{c1} holds by \eq{LBprob}. \hfill{$\qed$} \\

{\em Proof of \eq{voltail}.}
Taking $r$ so that $b_d r^d = t$ we have by \eq{0311e} that
$$
\BP [ V_d(C'(0,\Phi^0) ) > t ] \leq c_5 \exp (- c_5^{-1} r^{d/(d+1)} )
= c_5 \exp (- c_5^{-1} (t/b_d)^{1/(d+1)} )  
$$
while for the lower bound,
we have
\bean
\BP [ V_d (C'(0,\Phi^0)) \geq t ] \geq 
\BP [ \Phi (B_{2r}(0)) =0 ]
=
\exp (- 2^d t). \hfill{\qed}
\eean

\section{Central limit theorems}
\lbl{secCLT}
\eqco
In this section we derive  central limit theorems (CLTs) associated with
the lilypond model. 
We consider a sequence of binomial or Poisson point processes with finite
total number of points over an expanding sequence of bounded regions
in
 $\R^d$. Thus we consider, for example, 
 $ V_d(Z(\Phi_n))$ as defined in Section \ref{secintro},
  rather than
 $V_d(Z(\Phi) \cap  W_n)$. 
Our choice means that we can directly apply results in \cite{PY1} or
 \cite{PeEJP}, although
 it should be possible to 
 obtain similar results for $V_d(Z(\Phi) \cap W_n)$.

For $n,m \in \bN$,
recall the definition  of
 $X_1,X_2,\ldots$ and $\X_{n,m}$ at \eq{chidef}.
%
Set $W_n := n^{1/d} W$.
We consider the restricted  Poisson process 
$\Phi_n := \Phi \cap W_n$, and also
the binomial point process
$\X_n := \X_{n,n}$.
 For measurable $g: \R_+ \to \R$,
we give CLTs for sums of the form
$H_g(\Phi_n)$ and $H_g(\X_n)$, where 
for all $\varphi \in \bN$ with
$2 \leq \varphi(\R^d)< \infty$ we define 
$$
H_g(\varphi) := \sum_{x \in \varphi} g(\rho(x,\varphi)),
$$
and if $\varphi(\R^d)\in \{0,1\}$ we set $H_g(\varphi):=0$.
 For example if
$g(t) = b_d t^d$, then
by the hard-core property
$H_g(\Phi_n)$ is the total volume $V_d(Z(\Phi_n)){\bf 1}_{\{\Phi(W_n) \geq 2\}}$ and $H_g(\X_n)$ is 
$V_d(Z(\X_n)) $ for $n \geq 2$.

\begin{theorem}
\lbl{th:clt}
Suppose that there exists finite  $\beta >0 $ such that 
 $g:\R_+ \to \R$ satisfies the growth bound
\bea
\sup_{ t \in \R_+} 
 (1 + t^{\beta})^{-1}
|g(t)|
< \infty,
\lbl{ggrowth}
\eea
and that $g$ is not  Lebesgue-almost everywhere
constant.
Then there exist constants $0 < \tau_g^2 \leq \sigma_g^2 < \infty $ (dependent
on $g$ but independent
of the choice of the convex set $W$)  such that
as $n \to \infty$ we have
$n^{-1} \Var(H_g(\Phi_n)) \to \sigma_g^2$  and 
$n^{-1} \Var(H_g(\X_n)) \to \tau_g^2$,  and 
\bea
n^{-1/2} ( H_g(\Phi_n) - \BE H_g(\Phi_n)) \tod N(0,\sigma_g^2)  ; 
\lbl{CLT1}
\\
n^{-1/2} ( H_g(\X_n) - \BE H_g(\X_n)) \tod  N(0,\tau_g^2).
\lbl{CLT2}
\eea
\end{theorem}
The next theorem provides a similar CLT for the number of components.
Recall that $\kappa(A)$ (for $A \subset \R^d$)
 denotes the number of connected components
of $A$; 
for $\varphi \in \bN$ with $2 \leq \varphi(\R^d ) < \infty $  we define
$$
H_{\kappa}(\varphi) := \kappa \left( Z(\varphi) \right)
$$
which is also the number of  components of the graph $\G^*(\varphi)$,
and if $\varphi(\R^d)= 1$ we set $H_\kappa(\varphi):=1$,
and if $\varphi(\R^d)= 0$ we set $H_\kappa(\varphi):=0$.
\begin{theorem}
\lbl{th:cptclt}
There are constants $0 < \tau_\kappa^2 \leq \sigma_\kappa^2 < \infty$ such
that
$n^{-1} \Var(H_\kappa(\Phi_n)) \to \sigma_\kappa^2$  and 
$n^{-1} \Var(H_\kappa(\X_n)) \to \tau_\kappa^2$,  and 
\bea
n^{-1/2} ( H_\kappa(\Phi_n) - \BE H_\kappa(\Phi_n)) \tod N(0,\sigma_\kappa^2)  ; 
\lbl{CLT3}
\\
n^{-1/2} ( H_\kappa(\X_n) - \BE H_\kappa(\X_n)) \tod  N(0,\tau_\kappa^2).
\lbl{CLT4}
\eea
\end{theorem}

We shall prove Theorems \ref{th:clt} and \ref{th:cptclt} 
 by using  results from  \cite{PY1}
which have the advantage of showing that
the limiting variance is non-zero.
An alternative would be to use results from \cite{PeEJP}
 (an approach based instead on sub-exponential stabilization)
which could be used to  give a formula for $\sigma^2$
in terms of integrated two-point correlation functions,
and also to provide Gaussian limits for the random
measures associated with the sums $H_g(\Phi_n)$ and $H_g(\X_n)$;
 moreover, using Theorem 2.5 of \cite{PY2},
 it should be possible to provide Berry-Esseen type
error bounds associated with \eq{CLT3}   converging to zero at rate 
$O(n^{\varepsilon - 1/2})$ for any $\varepsilon >0$. 
However, we do not give details of these alternative approaches here.

The approach of \cite{PY1} is based on a notion of
{\em external} stabilization. Given a 
 real-valued functional $H(\varphi)$,
defined on finite $\varphi \in \bN$ in a 
translation-invariant manner (i.e. with $H(x+\varphi) = H(x)$
for all nonempty $\varphi \in \bN $ and all $x \in \R^d$), we say that
a nonnegative random variable $\tilde{R}$
is a {\em radius of external
stabilization} for $H$ if 
there is a further random variable $\tilde{\Delta}$ such that
$$
H((\Phi^0 \cap B_{\tilde R}(0)) \cup \psi) -
H(( \Phi \cap B_{\tilde R}(0) ) \cup \psi )
=  \tilde{\Delta}
$$
for all $\psi \in \bN$ with
 $\psi(\R^d)< \infty$ and $\psi (B_{\tilde{R}}(0))
=0$. If $H$ has  a radius of external stabilization
 $\tilde{R}$ with $\BP(\tilde{R} < \infty)=1$, we
say $H$ is {\em externally stabilizing} and refer
to $\tilde{\Delta}$ as the {\em add-one cost} of $H$.

To prove Theorems \ref{th:clt} and \ref{th:cptclt}
we shall use the following result, which is Theorem
2.1 of \cite{PY1}.
\begin{theorem} \lbl{PYthm} 
Suppose $H$ is externally stabilizing with add-one cost
$\tilde{\Delta}$, and satisfies
the moments condition
\bea
\sup_{n,m \in \N: n \geq 4, m \in [n/2,3n/2]} \sup_{x \in W_n} 
\BE ( | H(\X^x_{n,m}) - H( \X_{n,m}) |^4  ) < \infty
\lbl{umomH}
\eea
along with (for some $\beta \in(0,\infty)$) the growth bound
\bea
|H(\varphi) | \leq \beta (\diam(\varphi) + \varphi(\R^d))^\beta, ~~~ \varphi
\in \bN, ~ 0 < \varphi(\R^d) <  \infty. 
\lbl{Hgrowth}
\eea
Then there are constants $0 \leq \tau^2 \leq \sigma^2 < \infty$ such
that
$n^{-1} \Var(H(\Phi_n)) \to \sigma^2$  and
$n^{-1} \Var(H(\X_n)) \to \tau^2$,  and
\bean
n^{-1/2} ( H(\Phi_n) - \BE H(\Phi_n) ) \tod N(0,\sigma^2)  ;
\\
n^{-1/2} ( H(\X_n) - \BE H(\X_n) ) \tod  N(0,\tau^2).
\eean
Moreover,
if $\tilde{\Delta}$ has a non-degenerate distribution then $\tau^2 >0$
\end{theorem}


{\em Proof of Theorem  \ref{th:clt}.}
 We shall use Theorem \ref{PYthm}.  
Recall the definition \ref{0605a} of $\Rex(\varphi)$. For $\varphi \in \bN$ with
 $\Rex(\varphi) = 9r < \infty$, set
$$
\Delta_g(\varphi) := \left( 
\sum_{x \in \varphi^0 \cap B_{2r}(0) } g(\rho(x, \varphi^0))
\right)
- \sum_{x \in \varphi \cap B_{2r}(0) } g(\rho(x, \varphi)).
$$
By Lemma \ref{Rexlem}, if $\varphi \in \bN$ with $\Rex(\varphi) < \infty$ and
 $\psi \in \bN$ with $\psi(\R^d) < \infty$ and
$\psi ( B_{\Rex(\varphi)}(0) ) =0$, then
$$
H_g((\varphi^0 \cap B_{\Rex(\varphi)}(0) ) \cup  \psi) 
- H_g(\varphi \cap B_{\Rex(\varphi)}(0)  ) = \Delta_g(\varphi). 
$$
Hence, since $\Rex(\Phi) < \infty $ almost surely by
Lemma \ref{Rexfinlem}, $H_g$ is externally stabilizing.
Using \eq{ggrowth}, it is easy to see that  $H=H_g$ 
also satisfies the growth bound \eq{Hgrowth}.
We need to check that $H_g$ satisfies the moments
condition  \eq{umomH}, and that
the add-one cost $\Delta_g(\Phi)$ is non-degenerate.
We demonstrate these in Lemmas \ref{lmom} and \ref{lemnondeg} below.
Given these, we can apply Theorem  \ref{PYthm}
to get the result.
\hfill{$\qed$} \\

\begin{lemma}
\lbl{lmom}
Suppose $g: \R_+ \to \R$ satisfies the growth bound \eq{ggrowth}.
Then
the functional $H= H_g$ satisfies the moments condition \eq{umomH}.
\end{lemma}
{\em Proof.}
Since
$\chi_{n,m} \subset W_n$,
by Lemma \ref{Rexlem} (ii),
and Lemma \ref{Rexfinlem},
 for all  $r >0$, $n \geq 4$,   
$x \in W_n$ and $m \in [n/2,3n/2]$ we have 
\bea
\BP \left[  H_g(\chi_{n,m}^x)  - H_g(\chi_{n,m})  \neq
\left(
\sum_{z \in \chi_{n,m}^x \cap B_{2r}(x)} g(\rho(z,\varphi^x))
\right)
- \sum_{y \in \chi_{n,m} \cap B_{2r}(x)} g(\rho(y,\varphi))
 \right] \nonumber \\
 \leq c_4 \exp(- c_4^{-1} r^{d/(d+1)}).
\lbl{100118a}
\eea

Now observe that
 there is a constant $c_8$ such that 
by the assumed growth bound \eq{ggrowth} on $g$ and the bound
$\rho(x,\varphi) \leq D(x,\varphi)$,
 if $2 \leq \varphi(B_{2r}(x)) \leq 4 b_d (2r)^d$, then 
\bea
\left|
\left(
 \sum_{z \in \varphi^x \cap B_{2r}(x)} g(\rho(z,\varphi^x))
\right)
- \sum_{y \in \varphi \cap B_{2r}(x)} g(\rho(y,\varphi)) \right|
\leq c_8(1+ r^{d+\beta}).
\lbl{0108a}
\eea 
By \eq{convlow}, if $n \geq 4$ and $m \in [n/2,3n/2]$ and
$x \in W_n$ and $0 < r \leq n^{1/d}\diam(W)$, then 
\bea
\BE  \chi_{n,m}(B_{2r}(x)) = (m/n) V_d(B_{2r}(x) \cap W_n) \in
\left[ (2c_0)^{-1} (2r)^d, 2 b_d (2r)^d \right].
\label{0708a}
\eea
In particular, if 
 $(8c_0)^{1/d} \leq r \leq n^{1/d} \diam(W)$, then
$ \BE  \chi_{n,m}(B_{2r}(x)) \geq 2$.
 By a standard Chernoff-type tail estimate for  the binomial distribution
(see e.g. Lemma 1.1 of \cite{Pe}),
and \eq{0708a},
there is a constant $c_9 \in (0,\infty)$ such that if
 $(8c_0)^{1/d} \leq r \leq n^{1/d} \diam(W)$, then
\bea
1 - \BP[ 2 \leq \chi_{n,m} (B_{2r}(x)) \leq 4 b_d (2r)^d ] \leq
 2 \exp(-c_9^{-1} r^{d}).
\lbl{0108b}
\eea
Moreover, if $r > n^{1/d} \diam(W)$ then since $V_d(W)=1$,
 the Bieberbach (isodiametric)
inequality (see e.g. \cite{Pe}) yields
$$
4 b_d (2r)^d \geq 2^{2d+2} n b_d (\diam (W) /2)^d \geq 2^{2d+2} n,
$$
so if $2 \leq m \leq 3n/2$ then trivially \eq{0108b}  holds in this 
case too.
Combining  \eq{100118a}, \eq{0108a} and \eq{0108b}
 gives us, for all $n \geq 4,$ $ m \in [n/2,3n/2],$  $ x \in W_n$ 
and $r \geq (8c_0)^{1/d}$,  the tail bound
\bean
\BP [ | H_g(\X_{n,m}^x)  - H_g(\X_{n,m}) |
> c_8(1+ r^{d+\beta}) ]
\leq c_4 \exp(-c_4^{-1} r^{d/(d+1)}) + 2 \exp(-c_9 r^d),
\eean
which suffices to give us the moments bound \eq{umomH} for $H_g$. \hfill{$\qed$} \\

\begin{lemma}
\lbl{lemnondeg}
The random variable $\Delta_g(\Phi)$ has a nondegenerate distribution.
\end{lemma}
{\em Proof.}
By assumption $g$ is not almost everywhere constant, so by
the Lebesgue density theorem (as stated in e.g. \cite{Pe})
 there exist numbers
$a >b > 0$ and $ t_1 >0, $  $  t_2 >0 $ and $\eps \in (0, \min(t_1,t_2)/10)$
 such that $\int_{t_1}^{t_1+ \eps} {\bf 1}\{ g(t) \geq a\} dt >0$, and
$\int_{t_2}^{t_2+ \eps} {\bf 1}\{ g(t) \leq b\} dt >0$.
Choose $r >  \max(t_1,t_2)$ such that $\BP(\Phi \in F(99r,r) ) >0$
(this is possible by \eq{0311d}).

Enumerate the points of $\Phi$ as $X_1,X_2,\ldots$ with 
$|X_1| < |X_2| < \cdots$.
Set $D :=|X_1-X_2|$ and $R_0 : = 
\inf\{|x|: x \in B_{D/2}(X_1) \}$. 
For $i=1,2$ define the event
$$
E_i :=  \{ |D| \leq \eps \} \cap \{ R_0 \in (t_i,t_i+\eps) \}
\cap
\{|X_3| > 99  r \}  \cap \{ \Phi \in F(99r,r) \}
$$

Then $\BP (E_1) >0$ and $\BP(E_2) >0$. 
Also,
 given $E_i$ occurs the value of $\rho(0,\Phi^0)$ is equal to
$R_0$ and $H_g(\Phi^0 \cap B_s(0)) - H_g(\Phi \cap B_s(0)) = g(R_0)$
for any $s > 100r$, so (provided also $\Rex(\Phi) < \infty$)
  the value of $\Delta_g(\Phi)$ is  equal to 
$g(R_0)$.
 Moreover,
 given $E_i$, the distribution
of $R_0$ is absolutely continuous  on $(t_i,t_{i}+\eps)$
with a strictly positive density.
Therefore $\BP[\Delta_g(\Phi)  \geq a|E_1]$ and
$\BP[\Delta_g(\Phi)   \leq b |E_2]$ are both strictly positive.
Hence $\BP ( \Delta_g(\Phi) \geq a) >0$ and
$\BP ( \Delta_g(\Phi) \leq b) >0$.
 Thus $\Delta_g(\Phi)$ has a nondegenerate distribution
as required.  \hfill{$\qed$} \\

We now proceed towards a proof of Theorem \ref{th:cptclt}.
We wish to show the the functional $H_\kappa$ has an almost surely finite
 radius of external stabilization.
In fact we shall again use $\Rex$, as defined at \eq{0605a}.

Suppose $\Rex(\varphi) = 9r < \infty$.
Let $\varphi^* := \cup_{x \in \varphi \cap B_{3r}(0)} C(x,\varphi)$.
Let $N(\varphi)$ be the number of components of the subgraph of
$\G(\varphi) $ induced by the set of vertices
 $\varphi^*$. 
Let $N_0(\varphi)$
 be the number of components of the subgraph of
$\G(\varphi^0) $ induced by
 $\varphi^* \cup \{0\} $.  Let
$$
\Delta_{\kappa}(\varphi) :=  N_0(\varphi) - N(\varphi).
$$
The following lemma says that $\Rex$ 
 serves as a
radius of external stabilization for the functional $H_\kappa$:
\begin{lemma}
\lbl{lemserves}
(i) If $\varphi \in \bN$ with $\Rex(\varphi) < \infty$ and
 $\psi \in \bN$ with $\psi(\R^d) < \infty$ and
 $\psi (B_{\Rex(\varphi)}(0))=0$,
then
\bea
H_\kappa((\varphi^0 \cap B_{\Rex(\varphi)}(0) ) \cup  \psi) 
- H_\kappa((\varphi \cap B_{\Rex(\varphi)}(0)) \cup \psi  ) = \Delta_\kappa(\varphi). 
\lbl{0116b}
\eea
(ii) Let $n \in \N$.  Suppose $x \in W_n$, and
 $\varphi \in G_{n,r}(x).$ Then
$$
H_\kappa(\varphi^x)  - H_\kappa(\varphi) =
H_\kappa(\varphi^x \cap B_{7r}(x))  - H_\kappa(\varphi \cap B_{7r}(x)).
$$
\end{lemma}
{\em Proof.} Suppose $\Rex (\varphi) = 9r < \infty$.
 Then by  Lemma \ref{Rexlem}, the 
  radii of lilypond grains centred outside $B_{2r}(0)$
are unchanged when a point is inserted at 0, while
the radii of grains centred inside $B_{7r}(0)$
are unaffected by changes to $\varphi$ outside
$B_{9r}(0)$.

 Set $\varphi_1 := (\varphi \cap B_{9r}(0)) \cup \psi$. 
Let $x \in \varphi^* $ and 
$y \in \varphi_1 \setminus  \varphi^*$.
Then we claim that
 $\{x,y\} \notin \G^*(\varphi_1)$ and
$\{x,y\} \notin \G^*(\varphi_1^0)$.
Indeed, by Lemma \ref{prefroglem}
 we have $x \in B_{5r}(0)$, so
for $y \notin B_{7r}(0)$, 
 the claim follows from Lemma \ref{tri2lem} 
 and the fact that 
$\varphi \in F(7r,r/2) \cap F(5r,r/2)$
(so also 
$\varphi^0 \in F(7r,r/2) \cap F(5r,r/2)$).
In the case 
 $y \in B_{7r}(0)$, 
the claim follows from the definition of $\varphi^*$,
along with the fact that $\rho(y,\varphi_1) = \rho(y,\varphi)$
and $\rho(x,\varphi_1) = \rho(x,\varphi)$.
Similarly, $\{0,y\} \notin \G^*(\varphi_1^0)$.

Therefore the components of $\G^*(\varphi_1)$ induced by all 
$y \in \varphi_1 \setminus \varphi^*$
 do not meet the components containing all  $x \in \varphi^*$
 and this still holds after the addition of 0.
So the contribution of these components to
 $H_\kappa(\varphi_1^0) - 
 H_\kappa(\varphi_1)$ is zero and \eq{0116b} follows. 

Part (ii) is proved by much the same argument as for Part (i),
using Lemma \ref{fencenlem} instead of Lemma \ref{fencelem}
and Lemma \ref{tri2lem} (ii) instead of Lemma \ref{tri2lem} (i).
\hfill{$\qed$} \\

{\em Proof of Theorem \ref{th:cptclt}.}
By Lemmas \ref{Rexfinlem} and \ref{lemserves}, $H_\kappa$ is externally
 stabilizing.  It is easy to see that  $H_\kappa$ satisfies the growth
bound \eq{Hgrowth} with $\beta =1$.
  We need to check that $H_\kappa$ satisfies the 
 moments condition \eq{umomH}, and that
the limiting add-one cost $\Delta_\kappa(\Phi)$ is non-degenerate.
We demonstrate these in Lemmas \ref{lmom2} and \ref{lemnondeg2} below.
Given these, we can apply Theorem \ref{PYthm} to get the result.
\hfill{$\qed$} \\

\begin{lemma}
\lbl{lmom2}
The functional $H= H_\kappa$ satisfies the moments condition
\eq{umomH}.
\end{lemma}
{\em Proof.}
By Lemmas \ref{lemserves} (ii) and \ref{Rexfinlem},
for $n,m \in \N$ with $n \geq 4$ and $n/2 \leq m \leq 3n/2$,
 and $x \in W_n$  and $ r >0$, 
\bean
\BP [ H_\kappa(\X^x_{n,m} ) - H_\kappa(\X_{n,m}) \neq 
H_\kappa(\X_{n,m}^x \cap B_{7r}(x) )  
- H_\kappa(\X_{n,m} \cap B_{7r} (x) ) ]
\leq c_4 \exp (- c_4^{-1} r^{d/(d+1)} ).
\eean

Therefore, for all $\ell \in \N$,
\bea
\BP ( 
 | H_\kappa(\X^x_{n,m}) - H_\kappa( \X_{n,m}) | > \ell )
\leq 
c \exp ( - c^{-1} r^{d/(d+1)} )
+ \BP ( \X_{n,m}^x(B_{7r}(x)) > \ell ).
\lbl{0305c}
\eea
Now take $r : =r(\ell) := (1/7) (\ell/3b_d)^{1/d}$.
Then $\BE \chi_{n,m}(B_{7r}) \leq b_d (3/2) (7r)^d = \ell/2$.
Hence the last probability in \eq{0305c}
decays exponentially in $\ell$ by standard Chernoff tail estimates
for the binomial distribution (see e.g. Lemma 1.1 of \cite{Pe}),
so that overall there is a constant $c$ such that
the right hand side of \eq{0305c}
is bounded by  $c \exp (-c^{-1}\ell^{1/(d+1)})$, 
for all $\ell \leq n$,
 independently of $n $ and $m$.
Moreover, the left side of \eq{0305c} is zero for $\ell > n$.
This gives the uniform bound \eq{umomH} on fourth moments for $H_\kappa$.

\begin{lemma}
\lbl{lemnondeg2}
The limiting add-one cost $\Delta_\kappa(\Phi)$ has a non-degenerate 
distribution.
\end{lemma}
{\em Proof.}
Let $e \in \R^d$ be a unit vector. 
Enumerate the points of $\Phi$ as $X_1,X_2,\ldots$ with 
$|X_1| < |X_2| < \cdots$.
Given $r\geq 1$,
define events $E_{0,r}$ and $E_{1,r}$ as follows:
$$
E_{0,r} :=
 \{ \Phi(B_1(5e)) =2 \}
\cap \{|X_3| > 99r \} \cap \{\Phi \in F(99r,r)\};
$$
$$
E_{1,r} := \{\Phi(B_1(0))=1 \} \cap \{ \Phi(B_1(5e)) =2 \}
\cap \{|X_4| > 99r\} \cap \{\Phi \in F(99r,r)\}.
$$
Then by  \eq{0311d}, for large enough $r$ we have
$\BP(E_{0,r}) >0$ and $\BP(E_{1,r})>0$. But if $E_{0,r}$ occurs
and $\Rex(\Phi) < \infty$
then $\Delta_\kappa (\Phi) = 0$, while if $E_{1,r}$ occurs
and $\Rex(\Phi) < \infty$
 then
$\Delta_\kappa (\Phi) =1$.  This gives us the result.\hfill{$\qed$}

\section{Percolation theory for the enhanced model}
\lbl{secfrog}
\eqco

One interpretation (in $d=2$) of the absence of percolation in the Poisson
lilypond model, is that a frog is unable to travel infinitely
far by a continuous path along the lily pads. For the sake
of greater realism, it is natural to ask what happens if the frog
is allowed to jump.
 More mathematically, for $\delta >0$ we consider the
{\em enhanced union set}
 $Z^\delta$, where for any $A \subset \R^d$ we set $A^\delta := \cup_{x \in A}
B_\delta(x)$ so that in particular
\begin{align}
Z^\delta :=\bigcup_{x\in \Phi} B_{\rho(x,\Phi) + \delta}(x).
\end{align}
We investigate the connectivity properties of $Z^\delta$
(the parameter  $\delta$ represents half the distance
which the frog is able to jump).
In particular, we are concerned with the probability that 
$Z^\delta$ has an unbounded connected component.
Given $\delta$,
the event that this occurs is invariant under translations
of the Poisson process $\Phi$, so by  the ergodic property
of this Poisson process (see Proposition 2.6 of \cite{MR}), or
alternatively by an
 argument using the Kolmogorov zero-one law,
 this probability is either zero or 1, and it
is one if and only if
$$
\BP ( E_\infty(\delta)) >0
$$
where $E_\infty(\delta) $ denotes the event that there is an
infinite  component of $Z^\delta $ containing the origin.
Accordingly, we define the {\em critical enhancement}
$\delta_c : = \delta_c(d)$ by
$$
\delta_c := \inf \{ \delta >0: \BP (E_\infty(\delta)) >0 \}.
$$ 
Our main result in this section
 says that if the range of our jumping frog is sufficiently
small then it is still unable to travel infinitely far.
\begin{theorem}
\label{frogthm}
If $d =1$ then $\delta_c = \infty$.
If $d \geq 2$, then $0 < \delta_c < \infty$.
\end{theorem}
For $d=1$,
it is easy to see that
$\delta_c = \infty $.
For $d \geq 2$, the fact that $\delta_c$ is finite is immediate from the 
basic fact in continuum percolation,
 that for sufficiently large $r$ the union of balls
$\cup_{x \in \Phi}B_r(x)$ percolates (see for example \cite{Grimm}
or  \cite{MR}).
Thus, to prove Theorem \ref{frogthm} we need only to consider
the case with $d \geq 2$ and show $\delta_c >0$.

Let $D_r(x,\varphi)$ be the minimum of all non-zero pairwise  distances
between lilypond balls centred in $B_{8r}(x)$. That is, let 
$$
D_r(x,\varphi) = \min \{ |u -v| - \rho(u,\varphi) - \rho(v,\varphi)
 : u, v \in \varphi \cap B_{8r}(x),
|u-v| > \rho(u,\varphi) + \rho(v,\varphi) \}
$$
and set $D_r(x,\varphi) = + \infty $ if there are no such pairs $(u,v)$.
 Note that $D_r(x,\varphi)$ is   strictly positive
since, if finite, it is the minimum of a finite set of strictly 
positive numbers.

\begin{lemma}
\lbl{measlem}
Given $x \in \R^d$ and $r>0$, the event
$$
\{\Phi \in G_r(x) \} \cap \{D_r(x,\Phi) > \delta \} 
$$
is measurable with respect to $\sigma(\Phi \cap B_{9r}(x))$.
\end{lemma}
{\em Proof.}
Let $\varphi \in \bN$.
By Lemma \ref{stopsetlem},
$S(y,\varphi) \subset B_r(y) $ if and only if
$S(y,\varphi \cap B_r(y)) \subset B_r(y)$.
Hence,
 $\varphi \in G_r(x)$
if and only if $\varphi \cap B_{9r}(x) \in G_r(x)$. If
this is the case, then by Lemmas \ref{stopsetlem}
and \ref{l2.1}, we also have $D_r(x,\varphi) = D_r(x,\varphi
\cap B_{9r}(x))$. Hence the displayed event is identical to
 the event
$$
\{(\Phi \cap B_{9r}(x)) \in G_r(x) \} \cap \{D_r(x,\Phi \cap B_{9r}(x))
 > \delta\} 
$$
which  
is measurable with respect to $\sigma(\Phi \cap B_{9r}(x))$.
\hfill{$\qed$} \\

 A key claim is the following:

\begin{lemma}
\lbl{froglem}
Let $ \delta \in (0,1/2)$ and let $r >2$.
Let $x \in \R^d$.  If there is a continuous path
 in $Z^\delta$ from $\R^d \setminus B_{7r}(x)$ to $B_r(x)$,
then the event
$\{ \Phi \notin G_r(x) \} \cup \{ D_r(x,\Phi) \leq 2 \delta \}$
occurs.
\end{lemma}
{\em Proof.} 
Suppose that
 $\{ \Phi \notin G_r(x) \} \cup \{ D_r(x,\Phi)
 \leq 2 \delta \}$ does not occur. 
Let $T$ denote the union  of components of $Z $ which
intersect with $B_{2r}(x)$.
Any component of $Z$ with 
all its Poisson points outside $B_{3r}(x)$ is contained in $\R^d \setminus
 B_{2r}(x)$, by Lemma \ref{trilem} because  $ \Phi \in F(x,3r,r/2)$.
Hence by Lemma  \ref{prefroglem},
 \bea
T \subset
\cup_{y \in \Phi \cap B_{3r}(x)} C'(y,\Phi) \subset B_{6r}(x).
\label{0527a}
\eea


 Consider all lilypond balls centred at Poisson
points outside $T$.  Those centred inside $B_{8r}(x)$   
are at distance  more than  $2 \delta$ from $T$, because
$D_r(x,\Phi) > 2 \delta$  by assumption.
 Those centred outside 
$B_{8r}(x)$ do not intersect $B_{7r}(x)$ because  of the assumption
that $ \Phi \in F(x,8r,r/2)$ and Lemma \ref{trilem},
and so are also distant more than $2 \delta$ from
$T $,  
 since
(\ref{0527a}) holds and
$r > 2\delta$. 

Thus, the set $T$ is distant more than 
$2 \delta $ from $Z \setminus T$,  and hence $T^\delta$
is disconnected from the rest of $Z^\delta$. Finally,
by definition of $T$,
$T^\delta \subset B_{7r}(x)$ and
$Z^\delta \setminus T^\delta$ is disjoint from $B_r(x)$, both 
because $\delta < r$, so that there is no continuous 
path in $Z^\delta$ from $B_r(x)$ to $\R^d \setminus B_{7r}(x)$.
\hfill{$\qed$} \\

{\em Proof of Theorem \ref{frogthm}.}
Recall the definitions of 
$R(\varphi)$ and  $G_r(x)$ at \eq{0611a} and 
\eq{Grdef}  respectively.
By Lemma \ref{Rexfinlem},
$ \BP ( \Phi \notin G_r(0)) \to 0 $
as $r \to \infty$.
Also, $D_r(x)$ is a  strictly positive random variable,
Hence, given $\eps >0$, we can choose
$r >2 $ large enough and $\delta \in (0,1/2) $ small enough such that 
\bea
\BP (\{ \Phi \notin G_r(0) \} \cup \{D_r(0) \leq  2 \delta\} ) < \eps.
\lbl{0311c}
\eea
Note that
$\BP ( \{ \Phi \notin G_r(x) \}
 \cup \{D_r(x) \leq  2 \delta\} ) $  is the same for all $x$.

Now divide $\R^d$ into boxes (cubes) of side $2r d^{-1/2}$,
labelled $Q_z, z \in \Z^d$, by setting
$$
Q_z = 2rd^{-1/2}z + [-rd^{-1/2}, rd^{-1/2}]^d . 
$$
Define the random field $(Y_z,z \in \Z^d)$ by
$$
Y_z := 1- {\bf 1}_{
\{\Phi \notin G_r(2rd^{-1/2}z) \} \cup \{D_r(2rd^{-1/2}z,\Phi) \leq \delta\} }.
$$ 
Since  $Q_z \subset B_r( 2rd^{-1/2} z)  $, by Lemma \ref{froglem} 
if there is a 
 continuous path in $Z^\delta$ from $Q_z$ to $ \R^d \setminus
B_{6r}(2rd^{-1/d} z)$,
then $Y_z = 0$ almost surely.

If there is an infinite component in $Z^\delta$, then 
there must be an unbounded  continuous path in $Z^\delta$, and
by taking successive boxes along the path, there is an
infinite sequence   
$(z_1,z_2,z_3,\ldots)$
of distinct  elements
 with each $z_i \in \Z^d$ and $\|z_i -  z_{i+1}\|_\infty \leq 1 $ for each $i$,
 such that $Y_{z_i}=0$  for all $i$.
 
Given $z \in \Z^d$, by Lemma \ref{measlem} the random variable $Y_z$ 
is measurable with respect to
$\sigma(\Phi \cap B_{9r}(2rd^{-1/2}z))$.
Hence the random field $(Y_z,z \in \Z^d)$ 
is independent of 
$Y_{z'}$ for all all sites $z'$ with $2rd^{-1/2}|z'-z| > 18 r $,
i.e. with $|z' -z| > 9 d^{1/2}$.
Thus $Y_z$ is independent of $Y_{z'}$ whenever the graph
distance between $z$ and $z'$ exceeds $9d$. In fact,
$(Y_z,z \in \Z^d)$
 is a $9d$-dependent 
random field in the sense of
\cite{Grimm}.

Let $p_c$ be the  critical probability for site
percolation on the lattice with  vertex set $\Z^d$ and
edges between each pair $(z,z')$ with 
$\|z-z'\|_\infty =1$; it is well known \cite{Grimm} that $p_c >0$.
By (\cite{Grimm}, Theorem (7.65)), and  \eq{0311c},
we can choose $r$ to be so large  and $\delta$ to be so small that
the random field 
$(Y_z,z \in \Z^d)$  stochastically dominates a random field $(Y'_z, z \in \Z^d)$
consisting of independent  Bernoulli random variables
with $\BP(Y'_z=0) = p_c/2$ for each $z$. 
Thus, with this choice of $\delta$
there is almost surely no infinite path through
the lattice of sites with $Y_z=0$, and hence 
no infinite component in $Z^\delta$;
hence $\delta_c >0$ as asserted.
\hfill{$\qed$}


\end{document}